%% file: article.tex
\documentclass{amsart}

\usepackage{subfigure}
\usepackage{xspace}
\usepackage{amsmath}
\usepackage{amssymb}
\usepackage{amscd}
\usepackage{graphicx}
\usepackage{hyperref}
\usepackage{bbm}
\usepackage{tikz}

\usetikzlibrary{arrows,snakes,backgrounds}

\def\R{{\cal R}}

\def\R{{\mathbb R}}
\def\N{{\mathbb N}}
\def\S{{\bar \N}^{K-2}}
\def\E{\mathbb{E}}
\def\Pr{\mathbb{P}}

\def\bp{{\it Proof.}\ }
\def\ep{\hfill $\Box$}

\newcommand{\leqst}{\ensuremath{\leq_\mathrm{st}}}
\newcommand{\geqst}{\ensuremath{\geq_\mathrm{st}}}

\newcommand{\ind}{\ensuremath{\mathbbm{1}}}
\newcommand{\vect}[3]{\ensuremath{#1_{#2:#3}}}
\newcommand{\quasi}{\ensuremath{(\tilde{N}_{2:L}^\alpha(t))}}
\newcommand{\quasiscaled}[1]{\ensuremath{\tilde{N}_{2:L,\beta}^\alpha(#1)}}
\newcommand{\hadam}{\ensuremath{\odot}}
\newcommand{\satphi}[2]{\ensuremath{\breve{\phi}^{#1}_{#2}}}
\newcommand{\satpsi}[2]{\ensuremath{\breve{\psi}^{#1}_{#2}}}
\newcommand{\satN}[2]{\ensuremath{\breve{N}^{#1}_{#2}}}
\newcommand{\satscaledN}[2]{\ensuremath{\breve{N}^{#1}_{#2,\beta}}}
\newcommand{\satX}[2]{\ensuremath{\breve{X}^{#1}_{#2}}}
\newcommand{\satpi}[2]{\ensuremath{\breve{\pi}^{#1}_{#2,\beta}}}

\newcommand{\diff}{\mathop{}\mathopen{}\mathrm{d}}

\newtheorem{theorem}{Theorem}
\newtheorem{corollary}{Corollary}
\newtheorem{lemma}{Lemma}

\newtheorem{remark}{Remark}
\newtheorem{proposition}{Proposition}
\newtheorem{defi}{Definition}

\title[Stability of networks without congestion control]{On the flow-level stability of data networks without congestion control: the case of linear networks and upstream trees}
\author{Mathieu Feuillet}
\address{INRIA Paris-Rocquencourt, Domaine de Voluceau, 78153 Le Chesnay, France}
\email{mathieu.feuillet@inria.fr}
\date{\today}

\begin{document}
\maketitle

\begin{abstract}
In this paper, flow models of networks without congestion control are
considered. Users generate data transfers according to some Poisson processes
and transmit corresponding packet at a fixed rate equal to their access rate
until the entire document is received at the destination; some erasure codes are used
to make the transmission robust to packet losses. We study the stability of the
stochastic process representing the number of active flows in two particular cases:
linear networks and upstream trees. For the case of linear networks, we notably use
fluid limits and an interesting phenomenon of ``time scale separation'' occurs. Bounds
on the stability region of linear networks are given. For the case of upstream trees,
underlying monotonic properties are used. Finally, the asymptotic stability of
those processes is analyzed when the access rate of the users decreases to 0.
An appropriate scaling is introduced and used to prove that the stability region
of those networks is asymptotically maximized.
\end{abstract}

\section{Introduction}
\label{section:intro}
\input{intro.tex}

\section{Flow-level model}
\label{section:model}
\input{model.tex}

\section{Stability conditions for a linear network}
\label{section:line}
\input{line.tex}

\section{Asymptotic stability}
\label{section:asymptotic}
\input{asymptotic.tex}



\appendix
\section*{Appendix}
\section{Proof of Proposition \ref{prop:fluid-limits}}
\label{appendix:proof1}
\input{proof1.tex}

\end{document}

%% file: intro.tex
Most Internet traffic derives from the transfer of stored documents, corresponding
to texts, music and movies. This traffic is elastic in the sense that the transmission
rate varies depending on current network congestion. Each document transfer is a flow of
packets whose rate is usually regulated by the Transmission Control Protocol (TCP).
The rate is allowed to increase linearly until a packet loss is observed.
The rate is then divided by two. These two mechanisms ensure that network bandwidth
is shared approximately fairly between the varying number of flows in progress
at any time. TCP also ensures that lost packets are retransmitted.

There is, however, no constraint on users to actually implement TCP algorithm
and a more aggressive, greedy congestion control would be individually more profitable.
Indeed, new, more aggressive versions of TCP are emerging to better exploit
increasing link capacity \cite{rfc3649,cheng-04,lisong-04}. Some have suggested
that the packet retransmission function of TCP could be
replaced by forward error correction in the form of source coding \cite{snoeren-05}.
It consists in the use of erasure codes in order to encode files
before sending them into the network. As a consequence, if the file has been
divided in $K$ packets, the destination is able to rebuild the data with
$(1+\varepsilon)K$ packets. A typical class of erasure codes which is adapted
to source coding are the \emph{Digital fountains} \cite{mitzenmacher-04}.
The use of source coding enables the relaxation of control requirements and
facilitates aggressive packet transmission policies.

If the congestion avoidance of TCP is no longer universally used, there is a
danger that the Internet may again experience the congestion collapse observed
in its early days \cite{rfc896}. Packets dropped downstream and therefore
retransmitted needlessly encumber upstream links, amplifying and spreading
congestion. The consequence may be that the rate realized by concurrent flows is
much less than the optimum (in terms of utility maximization) and may even tend
to zero. The capacity of a network can be defined as the demand, flow arrival
rate $\times$ mean flow size, that can be sustained. It depends on the way
bandwidth is shared. It is important to understand
what happens to capacity when the assumption that users compliantly implement TCP
is relaxed. This is our objective in the present paper.

In this paper, we  assume users are greedy and do not implement congestion control. They send
data at the greatest possible rate determined by an external constraint  that we
assimilate to an access rate.  Bandwidth sharing is determined by this behavior.
We consider a flow-level model, as decribed in \cite{massoulie-00}, where the bandwidth
allocation changes instantly as new flows arrive or existing flows terminate.
Tail dropping is interpreted in this model such that, at each link, the output
rate of flows is proportional to their input rate to the link in question.
This is
consistent with the assumption that packets are all equally likely to be dropped.
The drop rate is such that the overall output rate is bounded by the link capacity.
Flows are categorized into classes defined by their route  and their
access rate. We assume flows arrive according to a Poisson process and have an
exponentially distributed size. The stochastic
process describing the number of flows in each class is a Markov process. We say
that the network is stable when this process is recurrent positive and unstable
otherwise.

\begin{figure}[h]
  \centering
  \begin{minipage}{0.45\linewidth}
  \centering
  \includegraphics[width=0.5\linewidth]{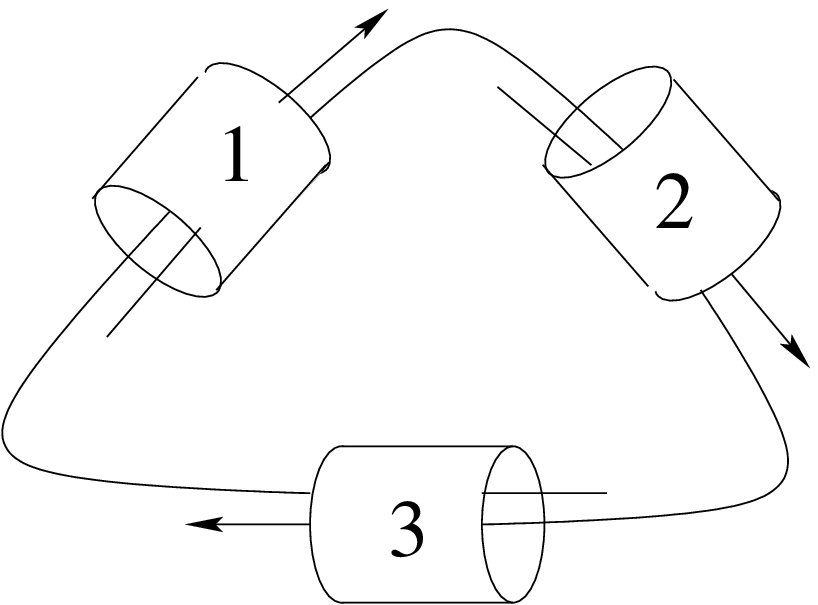}
  \caption{A cyclic network}
  \label{fig:cyclic}
  \end{minipage}
  \begin{minipage}{0.45\linewidth}
  \centering
  \newlength{\graphicheight}
  \settoheight\graphicheight{\includegraphics[width=0.5\linewidth]{triangle.eps}}
  \parbox[c][\graphicheight]{\linewidth}{\includegraphics[width=0.8\linewidth]{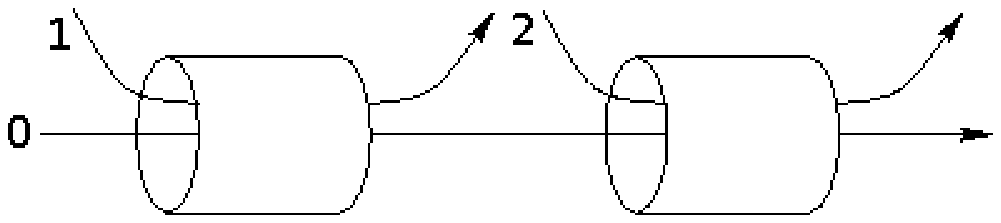}}\\
  ~\vfill
  \caption{A linear network}
  \label{fig:line}
  \label{fig:2-line}
  \end{minipage}
\end{figure}

To illustrate the realized bandwidth sharing without congestion control, consider
the linear network depicted in
Figure \ref{fig:line} with two links 1, 2 of capacity $C_1$ and $C_2$, respectively, $n_0$ class-0 flows
going through links 1 and 2, $n_1$ class-1 flows going through
link 1 and $n_2$ class-2 flows going through link 2. The access rate for class-$k$ flows is denoted by $a_k$.
The aggregate  input
rates of classes 0 and 1 at the first link are respectively $n_0a_0$ and $n_1a_1$ and
the total input rate at link $1$ is $R_1=n_0a_0+n_1a_1$. If $R_1>C_1$,
then the first link is saturated and the aggregate output rates of classes 0 and 1 are
respectively
$$
\theta_0^1 = C_1\frac{n_0a_0}{n_0a_0+n_1a_1}\quad \text{and}\quad
\theta_1^1 = C_1\frac{n_1a_1}{n_0a_0+n_1a_1}.
$$
If $R_1\leq C_1$ then the first link is not saturated and  $\theta_0^1=n_0a_0$
and $\theta_1^1=n_1a_1$. The aggregate input rate of class 0 at link 2 is $\theta_0^1$
and the aggregate input rate of class 2 is $n_2a_2$. As above, we derive
the respective aggregate output rates of class 0 and 2 after link 2:
$$
\theta_0^2=  \theta_0^1\min\left(1, \frac{C_2}{\theta_0^1+n_2a_2} \right)\quad \text{and}\quad
\theta_2^2 = n_2a_2 \min\left(1, \frac{C_2}{\theta_0^1+n_2a_2} \right).
$$

The question of the capacity of networks without congestion control has already
been addressed in \cite{bonald-09}. It was
shown that congestion collapse does occur in cyclic networks like that depicted in Figure
\ref{fig:cyclic} where the system is unstable for any positive demand.  In acyclic
networks, however, capacity is not reduced to 0 and this performance degradation is
significantly mitigated by the access rates. It was conjectured in \cite{bonald-09}
that capacity in acyclic networks increases as access rates tend to 0 and, in the
limit, is determined only by the optimal per-link stability conditions. In this
paper, we prove that the conjecture is indeed true for two particular acyclic
topologies, linear networks and upstream trees. For linear networks, we derive
bounds on stability conditions for arbitrary access rates.

\paragraph{Averaging phenomenon}
We evaluate stability conditions by applying fluid limit techniques. The end-to-end
class in a linear network plays a particular role. In the fluid limit, all other
flow classes are shown to empty in finite time for any initial conditions. When all these classes are empty,
a phenomenon of \emph{averaging} occurs and allow us to prove that the end-to-end class also empties in finite time.
In that situation, the rate at which this class evolves depends on the stochastic
evolution of the remaining local classes which in turn depends on the current
fluid state of the end-to-end class.  This dependence is decoupled by a time
scale separation argument. The rate at which the fluid component (the end-to-end flow)
evolves is infinitely slower than the rate at which the stochastic components
evolve. It is therefore possible to apply a \emph{quasi-stationary} model through
 which the \emph{average} impact of the local traffic can be evaluated.

This is a non-trivial example of a local equilibrium
of fluid limits. A simpler example of such a phenomenon can be found in
\cite[9.6, p271]{robert-03}. Averaging in fluid limits has also been considered
in the context of wireless networks (see \cite{bonald-04} and subsequent papers
\cite{simatos-10} and \cite{ganesh-10}). The main difference here is that the local equilibrium
depends on the state of the fluid limit. Note that this is somewhat related to the averaging
phenomenon considered for loss networks  \cite{kurtz-94}.

To evaluate the impact of access rates on capacity, we introduce an ad-hoc scaling.
Scaled processes are shown to converge to a simple deterministic process. Generally
speaking, the convergence of processes does not imply the convergence of
stationary distributions. However, in our case, we can prove this is true in two
specific cases exhibiting an interesting example of limit
inversion for stochastic processes. We use the same technique as in \cite{dumas-02}
to prove this inversion.
We use the convergence of stationary distributions combined with the averaging
phenomenon to prove the asymptotic optimality of linear networks. For upstream
trees, we use convergence and monotonicity.

The paper is structured as follows. In Section \ref{section:model}, we give a
complete description of the considered model. In Section
\ref{section:line}, we study the
stability conditions of linear networks with general access rates.
In Section \ref{section:asymptotic}, we study the stability of
linear networks and upstream trees when access rates decrease to $0$
and prove their capacity is asymptotically optimal.

%% file: model.tex
\subsection{Bandwidth allocation}

Consider a network of $L$ links. Denote by $C_l$ the capacity of link
$l$.  A number of flows compete to share the capacity of these links. The flows
are categorized into a set of $K$ classes. Each class-$k$ flow has an access
rate to the network denoted by $a_k$ and follows a route of length $d_k$
defined as an ordered set of distinct links $r_k=\{r_k(1),
r_k(2),\ldots,r_k(d_k)\}$. Let $x=(x_1,\ldots,x_K)$ be the vector of input
rates of all classes, $\psi_k(x)$ the aggregate throughput of
class-$k$ flows. We refer to the vector $\psi(x)=(\psi_1(x),\ldots,\psi_K(x))$
as the bandwidth allocation.

Here, we consider allocations that result when users are greedy and transmit
at their maximum input rate in the network without any congestion control.
There can be losses on each link of the network and we need to describe the
evolution of the rate of each flow going through the network.  Specifically,
we assume that class-$k$ flows initially transmit at full rate
$\theta_k^0(x) = x_k= n_ka_k$ where $n_k$ is the number of class-$k$ active flows
and we denote by $\theta_k^i(x)$ the class-$k$ rate at the {\it output} of the
$i$-th link on its route $r_k$, for
$i=1,2,\ldots,d_k$.   The actual throughput $\psi_k(x)$ of class $k$
corresponds to the aggregate throughput of class-$k$ flows at the
output of the last link on their route, $\psi_k(x)=\theta_k^{d_k}(x)$.
Due to potential loss at each link, we have
\begin{equation}
  \label{eq:evol}
  \psi_k(x)=\theta_k^{d_k}(x)\leq \theta_k^{d_k-1}(x)\leq \ldots\le
  \theta_k^1(x)
  \leq \theta_k^0(x)=x_k=n_ka_k.
\end{equation}
In particular, $\theta_k^{i-1}(x)-\theta_k^i(x)$ is the loss for
class $k$ occurring at the $i$-th link on its route.

The total rate at the output of each link cannot exceed the capacity of this
link, so that
\begin{equation}
  \label{eq:const2}
  \forall l,\ \sum_{k,i:r_k(i)=l} \theta_k^i(x)  \leq C_l.
\end{equation}

In order to characterize the allocation achieved in the absence of
congestion control, it remains to determine how output rates depend on
input rates at each link. This depends on the buffer management
policy.  In the considered fluid model, losses occur on saturated
links only when the total input rate exceeds the capacity; the total
output rate is then assumed to be equal to the capacity. In the following, we
consider a single buffer management policy, Tail Dropping, which simply consists
in dropping each incoming packet when the buffer is full. We deduce that, at the flow
level, the output rate of classes from the link is proportional  to the input
rate of classes to the link.

Specifically, let $R_l(x)$ be the total input rate of link $l$:
$$
R_l(x)=\sum_{k,i:r_k(i+1)=l} \theta_k^i(x).
$$
If $R_l(x)\leq C_l$ then link $l$ is not saturated and there are no losses;
the output rate is equal to the input rate for each class. If $R_l(x)>C_l$, then
the link is saturated and there are losses proportional to the input rate.
The output rate of any class $k$ then satisfies for all
$i=1,2,3,\ldots,d_k$,
\begin{equation}
  \label{eq:tail-drop}
  \theta_k^i(x) = \theta_k^{i-1}(x)
  \min\left({C_l\over R_l(x)},1\right)\quad\hbox{where } l=r_k(i).
\end{equation}

When there is no ambiguity
and when the access rates are fixed, the bandwidth allocation is determined
by the number of flows in each class. We denote by $n=(n_1,\ldots,n_K)$ the
vector of the numbers of flows in progress and by the vector $\phi(n)=(\phi_1(n),\ldots,\phi_K(n))$
the bandwidth allocation where, for each $k$, $\phi_k(n)=\psi_k(a\hadam n)$
where $\hadam$ denote the componentwise product.

For example, we can consider a linear network under tail dropping (see
Section \ref{section:line}) with two links of capacity $1$ with
access rates all equal to $1$ and $1$ flow in each class. We get
the following bandwidth allocation:
$$
\phi_0(1,1,1)=\frac{1}{3},\ \phi_1(1,1,1)=\frac{1}{2},\
\phi_2(1,1,1)=\frac{2}{3}.
$$

It is not obvious that (\ref{eq:tail-drop}) defines a unique allocation for all
networks. The definition is clearly non-ambiguous for {\it acyclic} networks,
i.e., if links can be numbered in such a way that each route consists of an
increasing sequence of link indexes. The allocation then directly
follows from applying (\ref{eq:tail-drop}) to all $l=1,\ldots,L$,
successively. For general networks, the allocation is still well-defined and unique
as proved in \cite{bonald-09}. In this paper, we consider only
acyclic networks.

\subsection{Flow dynamics}

We assume that class-$k$ flows are generated according to a Poisson
process of intensity $\lambda_k$ and have independent, exponentially
distributed sizes with mean $1/\mu_k$. We define the traffic intensity
of class-$k$ flows as $\rho_k=\lambda_k/\mu_k$.

Under the above assumptions, let $N(t)$ denote the number of flows at time $t$,
$(N(t),t\in\R_+)$ defines a Markov process with transition rates $\lambda_k$ from state
$n$ to state $n+e_k$ and $\phi_k(n)\mu_k$ from state $n$ to
state $n-e_k$, where $e_k$ denotes the vector with 1 in the $k$-th
component and 0 elsewhere.  We say that the network is {\it stable}
when this Markov process is ergodic. There is a well-known necessary stability condition
(see \cite{bonald-al-06}) which is
that the traffic intensity is less than the capacity for all links:
\begin{equation}\label{eq:optimal}
  \forall l, \quad \sum_{k:l\in r_k} \rho_k < C_l.
\end{equation}
An allocation will be called optimal if this condition is also sufficient.
In the rest of this paper, this condition will be further referred to as the
\emph{optimal stability condition}.

In the following, we study stability conditions for two specific classes
of networks. In Section \ref{section:line}, we consider linear networks.
In Section \ref{section:asymptotic}, we study the behavior of stability conditions
when the access rates are small for linear networks and upstream trees.

%% file: line.tex
\subsection{Linear networks}

We consider a linear network with $L$ links as depicted in Figure \ref{fig:2-line}.
For the sake of simplicity, we assume that the capacity of all links is
1 but all the results are true in the general case and can be obtained by adapting the notation.
The linear network has $K=L+1$ different classes of flows. Class-0 flows go through
all links. We label the links of the network from $1$ to $L$.
Link $l$ is the $l$-th link on the route of
class-0 flows. Class-$l$ flows go through link $l$ only.
Equation (\ref{eq:tail-drop}) can be rewritten as
follows for a given state $n$ and for $1\leq k \leq L$:
\begin{align*}
  \theta_0^k(x) &= \min \left(\theta_0^{k-1}(x), \frac{\theta_0^{k-1}(x)}{\theta_0^{k-1}(x)+x_k}\right)\ \mathrm{,}\\
  \theta_k^1(x) &= \min \left (x_k,
    \frac{x_k}{\theta_0^{k-1}(x)+x_k} \right )\mathrm{.}
\end{align*}

To simplify the presentation, we suppose that $a_0= 1$ and $a_1=1$.
In that case, a single flow of class 0 or 1 is enough
to saturate the first link and we have
$$
    \theta_0^1(n\hadam a) = \frac{n_0}{n_1+n_0}\quad\text{and}\quad
    \theta_1^1(n\hadam a) = \frac{n_1}{n_0+n_1}.
$$
All the results presented here remain true without this assumption. Only the proof
of Proposition \ref{prop:averaging} has to be adapted.

In the rest of this section, we study the stability conditions
of the stochastic process resulting from this linear network. The
optimal stability conditions \eqref{eq:optimal} are:
\begin{equation}\label{eq:line-optimal}
\rho_0+\rho_k<1\ \mathrm{for}\ 1\leq k \leq L.
\end{equation}

\subsection{General fluid limits}
\label{subsection:stability-2-K}

We show that classes $2,\dots,L$ are favored compared to
classes $0$ and $1$ in the sense that the ergodicity conditions
of classes $2,\dots,L$ do not depend on classes $0$ and $1$
but the contrary is not true.

In order to study the ergodicity of the system, we need to define fluid limits.
In the following, we consider the norm such that for $x\in\R^K$, $\|x\|=\sum_{i=1}^K |x_i|$.
Let $(m^i, i\in \N)$ be a sequence of $\N^K$ such that
$\lim_{i\rightarrow \infty} \|m^i\|=\infty$. Define for all $m \in
\N^K$, the process $(N^m(t))$ describing the system evolution under
tail dropping when starting at $N^m(0)=m$. A fluid limit $\bar{Z}$ of
the system is defined as an accumulation point of the laws of the processes in the set
$$
\mathcal{C} = \left\{(N^{m^i}(\|m^i\|t)/\|m^i\|,t\geq 0),i\in\N\right\}
$$
which is regarded as a subset of the space $D_{\R^N_+}([0,\infty))$ with
Skorohod topology (see \cite{billingsley-99}). It can be proved, as in
\cite[Prop 9.3 p 246]{robert-03},  that the set $\mathcal{C}$ is relatively compact
and that all the fluid limits $\bar{Z}$ are continuous.

In order to prove ergodicity of the process $(N(t))$, we just have to
prove that there exists a time $T$ such that for any initial condition of
$(\bar{Z}(t))$, for all $t\geq T$, $\bar{Z}(t)=0$ (see \cite{dai-95}).
On the contrary, in order to prove transience, we have to prove that there exists
a time $T$ such that after $T$, any fluid limit $(\bar{Z}(t))$ increases linearly
(see \cite{meyn-95}). For that purpose, we have to characterize the fluid limits.

First, we consider fluid limits such that there exists $2\leq k_0 \leq L$ with
$\bar{Z}_{k_0}(0)>0$.
At the stochastic level, for each $k\geq2$, class $k$ can use an arbitrarily large part
of link $k$ leaving almost nothing to class $0$; for any $\varepsilon>0$,
since $\theta_0^{k-1}(n \hadam a)\leq 1$, if $n_k$ is large enough, we have
$$
\phi_k(n) \geq 1-\varepsilon,\ \mathrm{and}\ \theta_0^k(n \hadam a) \leq \varepsilon.
$$
Thus, at fluid level, the throughput for $(\bar{Z}_k(t))$ is always $1$ if
$\bar{Z}_k(t)>0$. On the contrary, the throughput for $(\bar{Z}_0(t))$
is $0$ as long as there exists $k_0$ with $\bar{Z}_{k_0}(t)>0$.
We can prove, as in \cite[Prop 9.4, p247]{robert-03}, the following proposition.

\begin{proposition}
If $(\bar{Z}(t))$ is a fluid limit of the system and there exists some time
interval $[0,t_0]$ such that for any $t\in[0,t_0]$, there is some
$k_0$ in $\{2,\dots,L\}$ such that $\bar{Z}_{k_0}(t)>0$ then, almost surely,
$(\bar{Z}(t))$ is differentiable on $[0,t_0]$ and satisfies for all $t\in[0,t_0]$,
\begin{align}
  \dot{\bar{Z}}_0(t) &= \lambda_0,\label{eq:fluid-gen1}\\
  \dot{\bar{Z}}_1(t) &= \lambda_1 - \mu_1\frac{\bar{Z}_1(t)}{\bar{Z}_0(t)+\bar{Z}_1(t)},\label{eq:fluid-gen2}\\
  \dot{\bar{Z}}_k(t) &= (\lambda_k - \mu_k) \ind_{\{ \bar{Z}_k(t)>0\}}\
  \mathrm{for}\  2\leq k\leq L.\label{eq:fluid-gen3}
\end{align}
\label{prop:general-fluid-limits}
\end{proposition}

For $2\leq k \leq L$, $\bar{Z}_k(t)$ thus decreases linearly to $0$
if $\rho_k < 1$. Note also that $(\bar{Z}_k(t))$ for $2\leq k \leq L$
does not depend on classes $0$ and $1$. This illustrates the strong asymmetry
between classes $0$ and $1$ on one hand and classes $2,\dots,L$ on the other hand.

If there exists $2 \leq k_0 \leq L$ such that $\rho_{k_0}>1$ then it is clear
that the fluid limit will increase linearly to infinity and the process $(N(t))$
is transient. On the contrary, if $\rho_{k_0}<1$ for all $2\leq k \leq L$ then
there exists a time $T$ such that for $t\leq T$ and $2\leq k \leq L$,
$\bar{Z}_k(t)=0$. In order to study ergodicity and transience of $(N(t))$
we then have to characterize the fluid limits such that $\bar{Z}_k(0)=0$
for $2 \leq k \leq L$. This is the purpose of the next subsection.

\subsection{Quasi-stationary fluid limits}
\label{subsection:frozen}

Before studying fluid limits, we need to define the average
throughput of class $0$ in the quasi-stationary case. This corresponds to the
case where the number of flows in classes $0$ and $1$ is fixed and the number of flows
in classes $2,\dots,L$ varies. In the following, we denote by
$\vect{x}{2}{L}=(x_2,\dots,x_L)$ the input rates of classes $2,\dots,L$.
Fixing the number of flows in classes $0$ and $1$ is equivalent to considering
a fixed rate $\alpha$ for class $0$ after link $1$. Since link $1$ is
of capacity $1$, $\alpha \in[0,1]$. Define
$\tilde{\theta}_0^1(\alpha,\vect{x}{2}{L})=\alpha$ and, for $2\leq k \leq L$,
\begin{align}
  \tilde{\theta}^k_0(\alpha,\vect{x}{2}{L}) &= \min\left(
    \tilde{\theta}^{k-1}_0(\alpha,\vect{x}{2}{L}),
    \frac{\tilde{\theta}^{k-1}_0(\alpha,\vect{x}{2}{L})}
    {\tilde{\theta}^{k-1}_0(\alpha,\vect{x}{2}{L})+x_k} \right),\label{eq:thetatilde}\\
  \tilde{\psi}_k(\alpha,\vect{x}{2}{L}) &= \min\left( x_k,
    \frac{x_k}{\tilde{\theta}_0^{k-1}(\alpha,\vect{x}{2}{L})+x_k}
  \right).\label{eq:phitilde}
\end{align}
The quantity $\tilde{\theta}^k_0(\alpha,\vect{x}{2}{L})$ is the output rate of
class $0$ after link $k$ and the throughput of class $0$ is
$\tilde{\psi}_0(\alpha,\vect{x}{2}{L})=\tilde{\theta}^L_0(\alpha,\vect{x}{2}{L})$.
Moreover, $\tilde{\psi}_k(\alpha,\vect{x}{2}{L})$ is the throughput of class
$k$. For $0\leq k \leq L$ and for $\vect{n}{2}{L}\in \N^{K-2}$, we define
$\tilde{\phi}_k(\alpha, \vect{n}{2}{L}) = \tilde{\psi}_k(\alpha,\vect{n}{2}{L} \hadam a)$.

Denote by $\quasi$ the Markov process describing
the evolution of classes $2,\dots,L$ in the quasi-stationary case. It is
defined on $\N^{K-2}$ and its transition rates are, for $2\leq k\leq L$,
\begin{align*}
  n_k \rightarrow n_k+1 &: \lambda_k,\\
  n_k \rightarrow n_k-1 &: \mu_k \tilde{\phi}_k(\alpha,\vect{n}{2}{L}).
\end{align*}

\begin{proposition}
  Consider a linear network with $L$ links of capacity $1$ in the quasi-stationary
  case, i.e., when the number of flows in classes $0$ and $1$ is fixed and the
rate of class 0 after link 1 is $\alpha$.

  The Markov process $\quasi$ describing the evolution of classes $2,\dots,L$ is
  ergodic if, for $2\leq k \leq L$, $\rho_k<1$. It is transient if there exists $k_0$
  in $\{2,\dots,L\}$ such that $\rho_{k_0}>1$.
  \label{prop:ergo-2L}
\end{proposition}
\bp It is enough to see that for $2 \leq k \leq L$,
$$
\min\left(n_ka_k ,C_k\right) \geq \tilde{\phi}_k(\alpha,\vect{n}{2}{L}) \geq
\frac{n_ka_k}{C_{k-1}+n_ka_k}
$$.

Since $C_k=1$, if $\rho_k<1$ for all $k$, there exists $\varepsilon >0$ and
$\eta$ in $\N$ such that for all
$\vect{n}{2}{L}$ in $\N^{K-2}\setminus \{0,\dots,\eta\}^{K-2}$ and for $2\leq k\leq L$
such that $n_k \geq \eta$,
$$
\lambda_{k} - \mu_{k} \tilde{\phi}_{k}(\alpha,\vect{n}{2}{L}) \leq - \varepsilon
$$
and $\quasi$ is therefore ergodic.

Conversely, if there exists $k_0$ such that $\rho_{k_0}>1$, then there exists
$\varepsilon >0$ and $\eta$ such that for all
$\vect{n}{2}{L}$ with $n_{k_0} > \eta$,
$$
\lambda_{k_0} - \mu_{k_0} \tilde{\phi}_{k_0}(\alpha,\vect{n}{2}{L}) \geq \varepsilon
$$
and $\quasi$ is therefore transient.
\ep

When $\quasi$ is ergodic, denote its unique stationary distribution by
$\tilde{\pi}^\alpha$. The average throughput of class $0$ in the
quasi-stationary case is then defined as follows:
\begin{equation}
  \forall \alpha \in [0,1]\ \bar{\phi}_0(\alpha) =
  \E_{\tilde{\pi}^\alpha}(\tilde{\phi}_0(\alpha,.)) = \sum_{\vect{n}{2}{L}\in \N^{K-2}}
  \tilde{\pi}^\alpha(\vect{n}{2}{L}) \tilde{\phi}_0(\alpha,\vect{n}{2}{L}).
  \label{eq:phi0}
\end{equation}

The average $\bar{\phi}_0$ depends on the traffic intensities and access rates of classes $2,\dots,L$.
In order to establish ergodicity and transience conditions for the linear network
in Theorems \ref{theo:stability} and \ref{theo:unstability} below, we
need the continuity of the function $\bar{\phi}_0$ with respect to $\alpha$.

Subsequent developments rely on the following notion of stochastic domination.
\begin{defi}
Let $X$ and $Y$ be two random variables in a partially ordered
measurable space. We denote $X\leqst Y$ and say that $Y$ \textit{stochastically dominates}
$X$ if, for any positive non-decreasing measurable function $f$, we have
$\E[f(X)]\leq \E[f(Y)]$.

On $\R^n$, we use the usual coordinate-wise partial orders. If $x$ and $y$ are
in $\R^n$, $x\leq y$ if $x_k \leq y_k$ for all $k$. On the space $D(\R_+,\R^n)$,
we say that $x\leq y$ if $x_k(t) \leq y_k(t)$ for all $k$ and $t$.

For any two random variables $X$ and $Y$ with respective distributions $\pi_X$ and
$\pi_Y$ and such that $X \leqst Y$,  we say that $\pi_Y$ \textit{dominates}
$\pi_X$.
\end{defi}

In particular, if $(X(t))$ and $(Y(t))$ are two Markov processes on $\N^n$ such
that $(X(t)) \leqst (Y(t))$ and $(X(t))$ is irreducible, the ergodicity of
$(Y(t))$ implies the ergodicity
of $(X(t))$ and the transience of $(X(t))$ implies the transience of $(Y(t))$.
If $(Y(t))$ and $(X(t))$ are ergodic and $X(\infty)$ and $Y(\infty)$ are random
variables with the stationary distributions of $(X(t))$ and $(Y(t))$ as distributions
then $X(\infty) \leqst Y(\infty)$. For more details on stochastic domination,
see \cite{massey-87}.

\begin{lemma}
Consider a linear network with $L$ links of capacity 1 in the quasi-stationary case
and assume that $\rho_k<1$ for $2\leq k \leq L$, the function $\bar{\phi}_0$
defined by Equation \eqref{eq:phi0} is continuous with respect to $\alpha$ on $[0,1]$.
\end{lemma}
\bp First, we prove that $\alpha \mapsto \tilde{\phi}_k(\alpha,\vect{n}{2}{L})$ is
$1$-Lipschitz for all $\vect{n}{2}{L}$ and $k$.
By definition, $\alpha \mapsto \tilde{\theta}_0^1(\alpha,\vect{n}{2}{L})$ is 1-Lipschitz.
According to \eqref{eq:thetatilde} and \eqref{eq:phitilde}, $\alpha \mapsto \tilde{\theta}_0^2(\alpha,\vect{n}{2}{L})$
and $\alpha \mapsto \tilde{\phi}_2(\alpha,\vect{n}{2}{L})$ are defined as the minimum
of two $1$-Lipschitz functions and are then $1$-Lipschitz functions. The result
follows by recursion.

In particular, this is true for $\alpha \mapsto \tilde{\phi}_0(\alpha,\vect{n}{2}{L})$.
In order to prove the continuity of $\bar{\phi}_0$, we just have to prove that
if a sequence $(\alpha_i)$ of $[0,1]$ converges to $\alpha_\infty$,
then $\tilde{N}^{\alpha_i}_{2:L}(\infty)$ converges in distribution to
$\tilde{N}^{\alpha_\infty}_{2:L}(\infty)$ and these random variables have
$\tilde{\pi}^{\alpha_i}$ and $\tilde{\pi}^{\alpha_\infty}$ as respective distributions.

We first prove that $\{ \tilde{\pi}^\alpha, \alpha \in [0,1] \}$ is
tight. Note that for $2\leq k \leq L$ and $\alpha$
in $[0,1]$,
\begin{equation}
\tilde{\phi}_k(\alpha,\vect{n}{2}{L}) \geq \frac{n_ka_k}{1+n_ka_k}.
\label{eq:sto-dom}
\end{equation}

Define the process $(\hat{N}_{2:L}(t))$ with the following transition
rates for $2\leq k \leq L$:
\begin{align*}
  n_k \mapsto n_k+1 &: \lambda_k,\\
  n_k \mapsto n_k-1 &: \mu_k \frac{n_ka_k}{1+n_ka_k}.
\end{align*}
and such that $\hat{N}_{2:L}(0)=\tilde{N}^\alpha_{2:L}(0)$. The components of
$(\hat{N}_{2:L}(t))$ are independent. For $2\leq k \leq L$, we have $\rho_k<1$ and
there exist $\varepsilon>0$ and $\eta$ such that, for $n_k\geq \eta$, we have
$$
\lambda_k - \mu_k \frac{n_k a_k}{1+n_k a_k} < -\varepsilon,
$$
which implies that $(\hat{N}_{k}(t))$ is ergodic; the ergodicity of $(\hat{N}_{2:L}(t))$
follows. Using Equation \eqref{eq:sto-dom} and a standard coupling argument, see for instance \cite[Lemma 1]{borst-08}, we can
prove that, for $\alpha$ in $[0,1]$ and $2\leq k \leq L$,
$(\hat{N}_{k}(t))$ stochastically dominates $(\tilde{N}_k^\alpha(t))$. It implies
that $(\hat{N}_{2:L}(t))$ stochastically dominates $\quasi$ for any $\alpha \in [0,1]$.
This implies that the stationary
distribution $\hat{\pi}$ of $(\hat{N}_{2:L}(t))$ dominates any distribution
$\tilde{\pi}^\alpha$. In particular, we have
$$
\forall \alpha \in[0,1], \forall \kappa \in \R,\
\tilde{\pi}^\alpha(\N^{K-2}\backslash [0,\kappa]^{K-2}) \leq
\hat{\pi}(\N^{K-2}\backslash [0,\kappa]^{K-2}).
$$
Thus, the set $\{ \tilde{\pi}^\alpha, \alpha \in [0,1] \}$ is tight on $\N^{K-2}$
with the usual topology.

By tightness, we can suppose that the sequence
$(\tilde{\pi}^{\alpha_i})_{i\in\N}$ is convergent. We call
$\tilde{\pi}_\infty$ its limit. All we have to do now is to characterize
this limit and prove its uniqueness.

Let $f$ be a bounded function on $\N^{K-2}$. For any $\alpha$ in $[0,1]$,
let $\tilde{\Omega}^{\alpha}$ denote the infinitesimal generator of $\quasi$.
For all $\vect{n}{2}{L} \in \N^{K-2}$,
\begin{align*}
\tilde{\Omega}^\alpha(f)(\vect{n}{2}{L}) = &\sum_{k=2}^L
\lambda_k (f(\vect{n}{2}{L}+e_k) - f(\vect{n}{2}{L}))\\
&- \sum_{k=2}^L \mu_k \tilde{\phi}_k(\alpha,\vect{n}{2}{L}) (f(\vect{n}{2}{L}-e_k)-f(\vect{n}{2}{L})).
\end{align*}
Since  $\alpha \mapsto \tilde{\phi}_k(\alpha,\vect{n}{2}{L})$ is 1-Lipschitz for
all $\vect{n}{2}{L}$ and $k$ and $f$ is bounded, there  exists $\eta$ such that
\begin{equation}
  \forall i \in \N,\ \forall \vect{n}{2}{L} \in \N^{K-2},\
  |\tilde{\Omega}^{\alpha_\infty}(f)(\vect{n}{2}{L}) - \tilde{\Omega}^{\alpha_i}(f)(\vect{n}{2}{L})|
  \leq \eta |\alpha_\infty - \alpha_i|.
  \label{eq:control-gen}
\end{equation}

Because $(\tilde{\pi}^{\alpha_i})$ is tight, for any $\varepsilon$,
there exists $\kappa>0$ such that
\begin{equation}
  \forall i\in\N,\ \tilde{\pi}^{\alpha_i}([0,\kappa]^{K-2}) \geq 1 - \varepsilon.
  \label{eq:control-tight}
\end{equation}

Using equations \eqref{eq:control-gen}) and \eqref{eq:control-tight},
we have
$$
\lim_{i \rightarrow \infty} \int_{\N^{K-2}}
\tilde{\Omega}^{\alpha_i}(f)(\vect{n}{2}{L})\tilde{\pi}^{\alpha_i}(d\vect{n}{2}{L}) =
\int_{\N^{K-2}} \tilde{\Omega}^{\alpha_\infty}(f)(\vect{n}{2}{L})\tilde{\pi}_\infty(d\vect{n}{2}{L})
 = 0.
$$
We deduce that $\tilde{\pi}_\infty$ is an invariant distribution of
$(\tilde{N}^{\alpha_\infty}_{2:L}(t))$ and by uniqueness, we have
$\tilde{\pi}_\infty = \tilde{\pi}^{\alpha_\infty}$. Finally, we have
$$
\lim_{i \rightarrow \infty} \bar{\phi}_0(\alpha_i) =
\bar{\phi}_0(\alpha_\infty).
$$
\ep

We are now able to characterize the fluid limits $(\bar{Z}(t))$ such that $Z_{k}(0)=0$
for $2\leq k \leq L$. The next proposition is proved in the appendix.

\begin{proposition}
  Consider a linear network of $L$ links of capacity 1 and assume that $a_0=1$ and $a_1=1$.

  If $(\bar{Z}(t))$ is a fluid limit of the system such that $\bar{Z}_k(0)=0$ for $2\leq k \leq L$, then almost
  surely
  \begin{align}
    \dot{\bar{Z}}_0(t) &= \lambda_0 - \mu_0
    \bar{\phi}_0 \left(\frac{\bar{Z}_0(t)}{\bar{Z}_1(t)+\bar{Z}_0(t)}\right),\label{eq:fluid-1}\\
    \dot{\bar{Z}}_1(t) &= \lambda_1 - \mu_1
    \frac{\bar{Z}_1(t)}{\bar{Z}_1(t)+\bar{Z}_0(t)}\label{eq:fluid-2},\\
    \dot{\bar{Z}}_k(t) &= 0\quad \text{for } 2\leq k \leq L\label{eq:fluid-3}
  \end{align}
  hold for all $t\in\R_+$ where $\bar{\phi}_0$ is the average throughput of
  class $0$ in the quasi-stationary case defined by Equation \eqref{eq:phi0}.
  \label{prop:fluid-limits}
\end{proposition}

This proposition shows that, at the fluid time scale,
there is a separation of time scales between classes $0$ and $1$ and
$2,\dots,L$. From the point of view of classes $0$ and $1$, classes
$2,\dots,L$ are quasi-stationary and there is an averaging
because they evolve infinitely faster. From the point of view of classes
$2,\dots,L$, the ratio of classes $0$ and $1$ is constant. This is a further
illustration of the asymmetry between classes $0$, $1$ on the one hand and classes
$2,\dots,L$ on the other hand.

\subsection{Stability of a linear network}

Now that we have characterized the fluid limits of the system,
we can study the stability conditions of the linear network. Theorem \ref{theo:stability}
gives a sufficient condition for stability and Theorem \ref{theo:unstability}
gives a sufficient condition for transience of $(N(t))$.

\begin{theorem}
Consider a linear network of $L$ links of capacity $1$.\\
If $\rho_k <1$ for $2\leq k\leq L$ and
$$\rho_0 < \inf_{1-\rho_1\leq x \leq 1}  \bar{\phi}_0(x)$$
then the network is stable, \emph{i.e.},  $(N(t))$ is ergodic.
\label{theo:stability}
\end{theorem}
\bp In order to prove ergodicity, we just have to
prove that there exists a time $T$ such that for any initial condition of
$(\bar{Z}(t))$, for all $t\geq T$, $\bar{Z}(t)=0$ (see \cite{dai-95}).

We consider a general fluid limit $(\bar{Z}(t))$. As proved in Proposition
\ref{prop:general-fluid-limits}, as long as there exists $k_0$ such that
$\bar{Z}_{k_0}(t) >0$, $(\bar{Z}(t))$ satisfied Equations \eqref{eq:fluid-gen1},
\eqref{eq:fluid-gen2} and \eqref{eq:fluid-gen3}.
Since $\rho_k<1$ for $2\leq k \leq L$, there exists a finite time
$T_0$ such that, for all $k\geq 0$ and $t\geq T_0$, $\bar{Z}_k(t)=0$. According to the
strong Markov property, the study is reduced to the case where the initial
state of the fluid limit verifies $\bar{Z}_0(0)+\bar{Z}_1(0)=1$ and $\bar{Z}_k(0)=0$
for $2\leq k \leq L$. According to Proposition \ref{prop:fluid-limits}, $(\bar{Z}(t))$ satisfies \eqref{eq:fluid-1},
\eqref{eq:fluid-2} and \eqref{eq:fluid-3}.
In particular, for all $t\geq 0$ and $2\leq k \leq L$, we have $\bar{Z}_k(t)=0$,
we then just have to study the behavior of $(\bar{Z}_0(t),\bar{Z}_1(t))$.

We define the following function:
\begin{align*}
  \forall \alpha \in [0,1],\ f(\alpha) = \inf_{\alpha
    \leq x \leq 1} \bar{\phi}_0(x)
\end{align*}
This function is continuous, non-decreasing and satisfies
$f(\alpha) \leq \bar{\phi}_0(\alpha)$ for all $\alpha \in
[0,1]$.
We define $(F(t))$ such that $F_0(0) = \bar{Z}_0(0)$, $F_1(0)=\bar{Z}_1(0)$
and for all $t \geq 0$,
\begin{align}
  \dot{F}_0(t) &= \lambda_0 - \mu_0 f
  \left( \frac{F_0(t)}{F_0(t)+F_1(t)} \right)\label{eq:Zc0},\\
  \dot{F}_1(t) &= \lambda_1 - \mu_1
  \frac{F_1(t)}{F_0(t)+F_1(t)}.
  \label{eq:Zc1}
\end{align}

Since we have $\bar{\phi}_0(\alpha) \geq f(\alpha)$, we
can deduce for all $t\geq0$, $\bar{Z}_0(t) \leq F_0(t)$ and
$\bar{Z}_1(t) \leq F_1(t)$. This follows from the following.
\begin{itemize}
\item If $\bar{Z}_0(t) = F_0(t)$ and $\bar{Z}_1(t) =
  F_1(t)$ then $\dot{\bar{Z}}_0(t) \leq \dot{F}_0(t)$
  and $\dot{\bar{Z}}_1(t) \leq \dot{F}_1(t)$;

\item if $\bar{Z}_0(t) < F_0(t)$ and $\bar{Z}_1(t) =
  F_1(t)$ then $\dot{\bar{Z}}_1(t) \leq \dot{F}_1(t)$;

\item if $\bar{Z}_0(t) = F_0(t)$ and $\bar{Z}_1(t) <
  F_1(t)$ then $\dot{\bar{Z}}_0(t) \leq \dot{F}_0(t)$.
\end{itemize}

To prove the stability of $(\bar{Z}(t))$, it is enough to prove that
$(F(t))$ will return to $0$ in finite time.
Let $\alpha(t) = F_0(t)/(F_0(t)+F_1(t))$. We have
\begin{equation}
  \dot{\alpha}(t) =  \frac{1}{F_0(t)+F_1(t)}
  \bigl [ \left ( \lambda_0 - \mu_0 f(\alpha(t)) \right) (1- \alpha(t))
    - \left( \lambda_1 - \mu_1(1-\alpha(t)) \right) \alpha(t) \bigr].
  \label{eq:adot1}
\end{equation}
Let $\beta(t) = \left ( \lambda_0 - \mu_0 f(\alpha(t)) \right) (1- \alpha(t))
    - \left( \lambda_1 - \mu_1(1-\alpha(t)) \right) \alpha(t)$.
So that
\begin{equation}
\dot{\alpha}(t) =  \frac{\beta(t)}{F_0(t)+F_1(t)}.
  \label{eq:adot2}
\end{equation}
We also define  $ \alpha_0 = \max \{\alpha \in [0,1],\ f(\alpha) = \rho_0\}$.
Since $\rho_0<f(1-\rho_1)$, we have $\alpha_0 < 1-\rho_1$. Further, for
$t_0$ such that $\alpha(t)=\alpha_0$ and $t_1$ such that $\alpha(t_1)=1-\rho_1$,
\begin{align*}
\beta(t_0) &= \alpha_0(\mu_1(1-\alpha_0)-\lambda_1)>0,\\
\beta(t_1) &= \rho_1(\lambda_0 -\mu_0f(1-\rho_1))<0.
\end{align*}
Since $f$ is continuous and non-decreasing,
there exists $\kappa>0$ and $\eta>0$ such that if $\alpha(t)\leq \alpha_0+\eta$,
then $\beta(t) \geq \kappa$ and  if $\alpha(t)\geq 1-\rho_1-\eta$, then
$\beta(t) \leq -\kappa$.

As a consequence, we deduce from Equation \eqref{eq:adot2} that
\begin{align*}
\dot{\alpha}(t)\geq \frac{\kappa}{F_0(t)+F_1(t)}\ &\mathrm{if}\ \alpha(t)\leq \alpha_0+\eta,\\
\dot{\alpha}(t) \leq -\frac{\kappa}{F_0(t)+F_1(t)}\ &\mathrm{if}\ \alpha(t)\geq 1-\rho_1-\eta.\
\end{align*}
Thus, if $\alpha(0) \in [\alpha_0+\eta,1-\rho_1-\eta]$, $\alpha(t)\in
[\alpha_0+\eta,1-\rho_1-\eta]$ for all $t>0$.

For all $t\geq 0$, we have
\begin{align*}
F_0(t)+F_1(t) &\leq F_0(0)+F_1(0) + (\lambda_0+\lambda_1)t,\\
&\leq 1+ (\lambda_0+\lambda_1)t,
\end{align*}
and we deduce that if $\alpha(0)$ is in $[0,\alpha_0+\eta]$ and as
long as $\alpha(t)$ stays in $[0,\alpha_0+\eta]$, we have
\begin{align*}
  \alpha(t) &\geq \int_0^t \frac{\kappa ds}{1+(\lambda_0+\lambda_1)s}\\
  &\geq \frac{\kappa}{\lambda_0+\lambda_1}\log\left
    (1+(\lambda_0+\lambda_1)t \right).
\end{align*}
Similarly, if $\alpha(0)$ is in $[1-\rho_1-\eta,1]$ and $\alpha(t)$
stays in $[1-\rho_1-\eta,1]$, we have
\begin{align*}
  \alpha(t) &\leq 1- \int_0^t \frac{\kappa ds}{1+(\lambda_0+\lambda_1)s}\\
  &\leq 1- \frac{\kappa}{\lambda_0+\lambda_1}\log\left
    (1+(\lambda_0+\lambda_1)t \right).
\end{align*}
We then define
$$
T_1=\frac{1}{\lambda_0+\lambda_1}\exp\left(\frac{\max(\alpha_0+\eta,
    \rho_1+\eta)(\lambda_0+\lambda_1)} { \kappa } \right)
$$
and if $\alpha(0) \in [0,\alpha_0+\eta] \cup [1-\rho_1-\eta,1]$ then
$\alpha(T_1) \in [\alpha_0+\eta,1-\rho_1-\eta]$.  Finally, for all
$t\geq T_1$, $\alpha(t) \in [\alpha_0+\eta,1-\rho_1-\eta]$.

Now, we know that $\alpha(t)$ will reach $[\alpha_0+\eta,1-\rho_1-\eta]$ in finite
time and stay in that interval; we just have to study the behavior of $F$
when $\alpha(0)$ is in $[\alpha_0+\eta,1-\rho_1-\eta]$.
Using equations (\ref{eq:Zc0}) and (\ref{eq:Zc1}), we have that, if
$\alpha(0) \in [\alpha_0+\eta,1-\rho_1-\eta]$, then, for all $t\geq 0$,
\begin{align*}
  \dot{F}_0(t) &\leq - \mu_0 \left (
    f(\alpha_0+\eta) - \rho_0 \right ) < 0,\\
  \dot{F}_1(t) &\leq - \mu_1 \eta
\end{align*}
and the two components decrease at least linearly to $0$ in finite time.
There exists $T_2$ such that for any fluid limit, for all $t\geq T_2$,
$F(t)=0$ and thus $\bar{Z}(t)=0$.
\ep

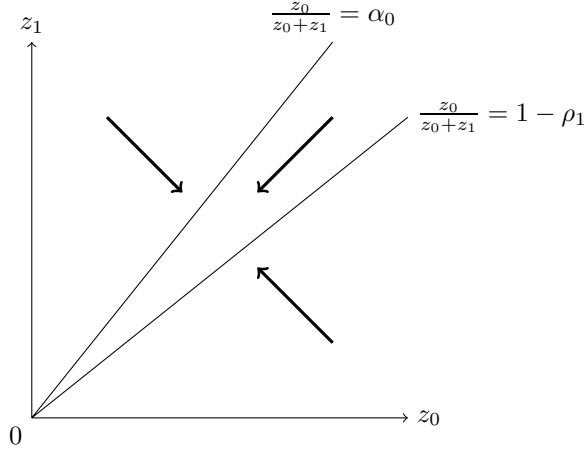
\begin{figure}
  \begin{center}
    \begin{tikzpicture}[scale=5]

      \node (0,0) [below left] {$0$};

      \draw[->,color=black] (0,0) -- (0,1) node [above] {$z_1$};
      \draw[->,color=black] (0,0) -- (1,0) node [right] {$z_0$};

      \draw (0,0) -- (0.8,1) node [above] {$\frac{z_0}{z_0+z_1} =
        \alpha_0$} (0,0) -- (1,0.8) node [right] {$\frac{z_0}{z_0+z_1}
        = 1-\rho_1$};

      \draw[very thick,->,color=black] (0.8,0.8) -- (0.6,0.6);
      \draw[very thick,->,color=black] (0.2,0.8) -- (0.4,0.6);
      \draw[very thick,->,color=black] (0.8,0.2) -- (0.6,0.4);

    \end{tikzpicture}
  \end{center}
  \caption{Dynamics of $(F_0(t),F_1(t))$}
  \label{fig:line-stable}
\end{figure}

The dynamics of $(F_0(t),F_1(t))$ described in this proof
is represented in Figure \ref{fig:line-stable}. Theorem \ref{theo:stability}
gives a sufficient condition for stability.  The next theorem gives a sufficient
condition for transience. The arguments of the proof are very similar.

\begin{theorem}
Consider a linear network of $L$ links of capacity $1$.
If there exists $k_0$ in $\{0,\dots,L\}$ such that $\rho_{k_0}>1$ or,
if for all $k$ in $\{0,\dots,L\}$, $\rho_k<1$ and
$$
\rho_0 > \sup_{0\leq x \leq 1-\rho_1} \bar{\phi}_0(x),
$$
then the network is unstable, \emph{i.e.}, $(N(t))$ is transient.
\label{theo:unstability}
\end{theorem}
\bp In order to prove the transience of the process, we just have to
prove that there exists a time $T$ such that after $T$, any fluid
limit $(\bar{Z}(t))$ increases linearly (see \cite{meyn-95}).

We consider a fluid limit $(\bar{Z}(t))$. If there is a $k_0$ such that
$\bar{Z}_{k_0}(0)>0$ then Equations (\ref{eq:fluid-gen1}), (\ref{eq:fluid-gen2})
and (\ref{eq:fluid-gen3}) are valid. It is obvious that, if there
exists $k_0$ such that $\rho_{k_0}>1$ then $(\bar{Z}_k(t))$ increases linearly
and the process $(N(t))$ is transient.

We suppose now that $\rho_k<1$ for all $k$ in $\{0,\dots,L\}$. In that case,
there is a time $T_0$ such that, for all $k$ in $\{2,\dots,L\}$ and $t\geq T_0$,
$\bar{Z}_k(t)=0$. According to the strong Markov property, we just have to study
the fluid limits such that $\bar{Z}_k(0)=0$ for $2\leq k \leq L$. In that case,
$\bar{Z}_{k}(t) = 0$ for $2\leq k \leq L$ and we just have to study the dynamics
of $(\bar{Z}_0(t), \bar{Z}_1(t))$.

We define the following function
$$
g(\alpha) = \sup_{0\leq x \leq \alpha} \bar{\phi}_0(x).
$$
This function is continuous, non-decreasing and, for any $\alpha$ in
$[0,1]$ satisfies $g(\alpha) \geq \bar{\phi}_0(\alpha)$. We also define the
process $(G(t))$ such that $G(0)=\bar{Z}(0)$ and
\begin{align}
  \dot{G}_0(t) &= \lambda_0 - \mu_0 g \left (
    \frac{G_0(t)}{G_0(t)+G_1(t)} \right ),
  \label{eq:Zh0}\\
  \dot{G}_1(t) &= \lambda_1 - \mu_1
  \frac{G_0(t)}{G_0(t)+G_1(t)}.
  \label{eq:Zh1}
\end{align}
As previously, we can prove that, for all $t\geq0$, $G_0(t) \leq \bar{Z}_0(t)$ and
$G_1(t) \leq \bar{Z}_1(t)$. So, we just have to prove that the
process $(G(t))$ increases linearly to infinity in order to prove the transience
of $(N(t))$.

We define $\alpha_0 = \min \{ \alpha \in [0,1],\ g(\alpha) = \rho_0\}$.
Since $\rho_0>g(1-\rho_1)$, we have $\alpha_0>1-\rho_1$. Define
$\alpha(t) = G_0(t)/(G_0(t)+G_1(t))$.

As in the proof of Theorem \ref{theo:stability}, we can show that there exists
$\eta>0$ such that, for any initial state $\alpha(0)$, after some finite time
$T_1$, $\alpha(t)$ reaches and stays in $[1-\rho_1+\eta,\alpha_0-\eta]$. We then
just have to study the dynamics of $G$ when $\alpha(0)$ is in
$[1-\rho_1+\eta,\alpha_0-\eta]$.

Using equations (\ref{eq:Zh0}) and (\ref{eq:Zh1}), we have that, if
$\alpha(t) \in [1-\rho_1+\eta,\alpha_0-\eta]$, then
\begin{align*}
  \dot{G}_0(t) &\geq  \mu_0 \left ( \rho_0 - g(\alpha_0-\eta)   \right ) > 0,\\
  \dot{G}_1(t) &\geq \mu_1 \eta.
\end{align*}
Thus, for any initial conditions, for all $t\geq T_1$, $G(t)$
increases linearly to infinity and the theorem is proved.
\ep

\begin{figure}
  \begin{center}
    \begin{tikzpicture}[scale=5]

      \node (0,0) [below left] {$0$};

      \draw[->,color=black] (0,0) -- (0,1) node [above] {$z_1$};
      \draw[->,color=black] (0,0) -- (1,0) node [right] {$z_0$};

      \draw (0,0) -- (0.8,1) node [above] {$\frac{z_0}{z_0+z_1} =
        1-\rho_1$} (0,0) -- (1,0.8) node [right] {$\frac{z_0}{z_0+z_1}
        = \alpha_0$};

      \draw[very thick,->,color=black] (0.6,0.6) -- (0.8,0.8);
      \draw[very thick,->,color=black] (0.2,0.8) -- (0.4,0.6);
      \draw[very thick,->,color=black] (0.8,0.2) -- (0.6,0.4);

    \end{tikzpicture}
  \end{center}
  \caption{Dynamics of $(G_0(t),G_1(t))$}
  \label{fig:line-unstable}
\end{figure}
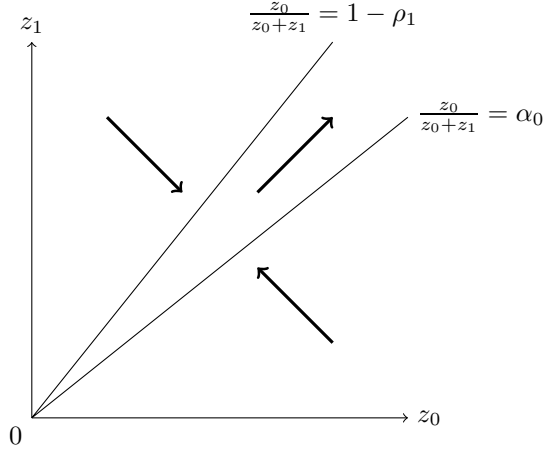

The dynamics of $(G_0(t),G_1(t))$ described in this proof
is represented in Figure \ref{fig:line-unstable}. In particular, we can remark
that when the network is unstable, there is always a time after which both
classes 0 and 1 of the fluid limit increase linearly to infinity.

\begin{remark}
  The stability region of the network is not reduced to $\{0\}$ since
$\inf_{1-\rho_1 \leq x \leq 1} \bar{\phi}_0(x) >0 $ when $\rho_1<1$ and it is
strictly included in the optimal stability region since
$\sup_{0\leq x \leq 1-\rho_1 } \bar{\phi}_0(x) < \min_k (1 - \rho_k)$.
\end{remark}

\begin{corollary}
  Consider a linear network with two links of capacity $1$.

  This network is stable if $\rho_k<1$ for all $0\leq k \leq 2$ and
  $\rho_0 < \bar{\phi}_0(1-\rho_1)$.

  This network is unstable if there exists $k_0$ in $\{0,1,2\}$ such that
  $\rho_{k_0}>1$ or if $\rho_k<1$ for $0\leq k\leq 2$ and
  $\rho_0 > \bar{\phi}_0(1-\rho_1)$.
\end{corollary}
\bp All that we have to prove is that the function $\bar{\phi}_0$ is
strictly increasing.

Consider the functions $\tilde{\phi}_0$, $\tilde{\phi}_2$, the processes
$(\tilde{N}^\alpha_2(t))$ and their stationary distribution $\tilde{\pi}^\alpha$
as defined in section \ref{subsection:frozen}. We assume that $\rho_2<1$.

The process $(\tilde{N}^\alpha_2(t))$ is one-dimensional and for $\alpha_1 \leq \alpha_2$,
 for all $n_2$, $\tilde{\phi}_k(\alpha_1,n_2) \geq \tilde{\phi}_k(\alpha_2,n_2)$.
This implies
$(\tilde{N}_2^{\alpha_2}(t)) \geq_{\mathrm{st}} (\tilde{N}_2^{\alpha_1}(t))$
and $\tilde{\pi}^{\alpha_2}$ dominates $\tilde{\pi}^{\alpha_1}$.

For any $\alpha\in[0,1]$ and any $n_2\in\N$,
$$
\tilde{\phi}_0(\alpha,n_2)+\tilde{\phi}_2(\alpha,n_2)=\min(n_2a_2+\alpha,1).
$$
The function $n_2 \mapsto \tilde{\phi}_0(\alpha,n_2)+\tilde{\phi}_2(\alpha,n_2)$
is therefore non-decreasing.

Applying the definition of stochastic domination, we have
$$
\E_{\tilde{\pi}^{\alpha_1}}(\tilde{\phi}_0(\alpha_1,.)+\tilde{\phi}_2(\alpha_1,.))
\leq
\E_{\tilde{\pi}^{\alpha_2}}(\tilde{\phi}_0(\alpha_1,.)+\tilde{\phi}_2(\alpha_1,.)).
$$

For any $n_2\in\N$, we have $\tilde{\phi}_0(\alpha_2,n_2)+\tilde{\phi}_2(\alpha_2,n_2)
\geq \tilde{\phi}_0(\alpha_1,n_2)+\tilde{\phi}_2(\alpha_1,n_2)$. Since
$\alpha_1<\alpha_2\leq1$, we have $\phi_0(\alpha_2,0)>\phi_0(\alpha_1,0)$ and we
can conclude
$$
\E_{\tilde{\pi}^{\alpha_1}}(\tilde{\phi}_0(\alpha_1,.)+\tilde{\phi}_2(\alpha_1,.)) <
\E_{\tilde{\pi}^{\alpha_2}}(\tilde{\phi}_0(\alpha_2,.)+\tilde{\phi}_2(\alpha_2,.)).
$$
Using the fact that
$\E_{\tilde{\pi}^{\alpha}}(\tilde{\phi}_2(\alpha,.))=\rho_2$, we can conclude
that if $\alpha_1<\alpha_2$, then $\bar{\phi}_0(\alpha_1) < \bar{\phi}_0(\alpha_2)$.
\ep

%% file: asymptotic.tex
From the results of the previous section, we know that ergodicity conditions in general
depend on the access rates $a_k$. In this section, we qualify this dependence
as the access rates become asymptotically small. In the next subsection, we introduce
a scaling on the access rates in order to understand the qualitative behavior of
networks when the access rates become small. In subsection \ref{subsection:line-asym},
we use this scaling to determine the limit of $\bar{\phi}_0$ when the access
rates tend to 0 and then deduce that linear networks are asymptotically
optimal. In subsection \ref{subsection:tree-asym}, we use the same scaling to
determine the asymptotic optimality of upstream trees. Subsections \ref{subsection:line-asym}
and \ref{subsection:tree-asym} are independent.

\subsection{Scaling over the access rates and flow sizes}
\label{subsection:scaling-access}

We fix all network parameters except the
arrival rates, the mean size of flows and the access rates. For each
$\beta\in \N$, define the process $(N_\beta(t))$ describing the
evolution of the number of flows in each class in a network where the access rate
of each class $k$ becomes $a_k/\beta$,
the arrival rate becomes $\lambda_k \beta$ and the inverse of the mean size
$\beta \mu_k$.

In this proof and in the proof in the appendix, we will need the definition of an
increasing process of a martingale. If $(M(t))$ is a square-integrable martingale on $\R$ null at 0,
$(\langle M \rangle_t)$ is the (essentially unique) increasing process such
that $(M^2(t)-\langle M \rangle_t)$ is a martingale. If $(M(t))$ and $(N(t))$ are
two square-integrable martingales on $\R$ null at 0, $(\langle M,N \rangle_t)$
is the (essentially unique) increasing process such that $(M(t)N(t) - \langle M,N \rangle_t)$
is a martingale. For more information on martingales and increasing processes,
see \cite{rogers-94}.

\begin{theorem}
  Consider an acyclic network with $L$ links and $K$ classes.

  If $$\lim_{\beta \to \infty} \frac{1}{\beta} a\hadam N_\beta(0) = \bar{X}(0),$$ then
  $(\beta^{-1} a \hadam N_\beta(t))$ converges in probability uniformly on compact sets
  to $(\bar{X}(t))$ such that
  \begin{equation}
    \dot{\bar{X}}_k(t) = a_k\lambda_kt - a_k\mu_k \psi_k(\bar{X}(t)),\quad \text{for }
    1\leq k \leq K.
    \label{limit-equation}
  \end{equation}
  \label{theo:large-numbers}
\end{theorem}
The convergence mentioned in this theorem is the uniform convergence in probability
on compact sets for the space $D_{\R^N_+}([0,\infty))$ with Skorohod topology
(see \cite{billingsley-99}).
\bp
We define the process $(M_\beta(t))$ such that
\begin{equation}
M_{\beta,k}(t) =a_k\left(\frac{N_{\beta,k}(t)}{\beta} - \frac{N_{\beta,k}(0)}{\beta}  - \lambda_k t + \mu_k
\int_0^t\psi_k\left(\frac{1}{\beta}a \hadam N_\beta(s)\right)\diff s\right),
\label{eq:scaling-sde-martingale}
\end{equation}
for $1\leq k \leq K$. The martingale characterization (see \cite{rogers-00}) of the Markov
jump process $(\beta^{-1} a \hadam N_\beta(t))$ shows that $(M_\beta(t))$ is a
martingale and its increasing processes  are given by, for $1\leq k,l \leq L$ and $k\neq l$,
\begin{align}
  \langle M_{\beta,k} \rangle_t &= \frac{a_k}{\beta} \left
    (\lambda_k t + \mu_k \int_0^t \psi_k\left(\frac{1}{\beta} a \hadam N_\beta(s)\right) \diff s \right ),
  \label{eq:quadratic-variation-1}\\
  \langle M_{\beta,k}, M_{\beta,j} \rangle_t &= 0.
  \label{eq:quadratic-variation-2}
\end{align}
Since $\psi_k$ is bounded by a constant $C$, we have
$\langle M_{\beta,k} \rangle_t \leq a_k (\lambda_k+C\mu_k)t/\beta$, for $t \geq 0$.

As $(M_\beta(t))$ is a martingale, we can use Doob's inequality. For any $t>0$
and $\varepsilon > 0$,
\begin{align*}
  \Pr \left ( \sup_{0\leq s \leq t} \|M_\beta\|_\infty(s) \geq
    \varepsilon \right )
  &\leq \sum_{k=1}^K \Pr \left ( \sup_{0\leq s \leq t} \|M_{\beta,k}\|_\infty(s) \geq \varepsilon \right )\\
  & \leq \frac{K}{\varepsilon^2} \max_k \E \left ( {M_{\beta,k}}^2(t) \right ),\\
  & \leq \frac{K}{\varepsilon^2} \max_k \langle M_{\beta,k} \rangle_t,\\
  & \leq \frac{Kt}{\varepsilon^2\beta}\max_k a_k (\lambda_k+C\mu_k).\\
\end{align*}
Thus, $(M_\beta(t))$ converges to $0$ in probability uniformly on
any compact when $\beta\rightarrow+\infty$.

We now consider the random variable $\sup_{0 \leq s \leq t} \left \|
  \beta^{-1} a \hadam N_\beta(t) - \bar{X}(t) \right \|_\infty$ and show that it
converges in distribution to $0$. Define
$Z_\beta(t) = \beta^{-1} a \hadam N_\beta(t) - \bar{X}(t)$. Using
(\ref{eq:scaling-sde-martingale}),
\begin{align*}
Z_{\beta,k}(t) = &\frac{1}{\beta} a \hadam N_{\beta,k}(0) -\bar{X}_k(0) + M_{\beta,k}(t)\\
&- \mu_k \int_0^t \left (\psi_k\left(\frac{1}{\beta} a \hadam N_\beta(s)\right) -
  \psi_k(\bar{X}(s)) \right ) \diff s,\quad \text{for}\ 1\leq k \leq K.
\end{align*}

The $\psi_k$ are all $\kappa$-Lipschitz for some $\kappa>0$ and, for $1\leq k\leq K$,
\begin{equation}
\begin{split}
  \sup_{0\leq s \leq t} |Z_{\beta,k}(s)| \leq
 & \left|\frac{a_k}{\beta} N_{\beta,k}(0) -\bar{X}_k(0)\right| +\sup_{0\leq s \leq t} |M_{\beta,k}(s)|\\
& + \kappa \mu_k \int_0^t
  \sup_{0 \leq u \leq s} |Z_\beta(u)|\diff s.
\end{split}
  \label{sup-zn}
\end{equation}

Define the function
$$
f_\beta(t) = \E \left ( \sup_{0\leq s \leq t} \|Z_\beta(s) \|_\infty
\right ).
$$

Using \eqref{sup-zn}, \eqref{eq:quadratic-variation-1} and
\eqref{eq:quadratic-variation-2}, we deduce the inequality for all
$s<t$,
$$
f_\beta(s) \leq \left\|\frac{1}{\beta}a \hadam N_\beta(0) -\bar{X}(0)\right\|_\infty + \frac{\hat{a}}{\beta} \left
  (\hat{\lambda}+ \hat{\mu} \right ) s + \hat{\mu} \int_0^s
f_\beta(u)\, \diff u
$$
where $\hat{a} = \max_k a_k$, $\hat{\lambda} = \max_k \lambda_k$ and
$\hat{\mu} = \max_k \mu_k$.

Applying Gronwall's lemma (see \cite{gronwall-19}),
$$
f_\beta(t) \leq \left ( \left\|\frac{1}{\beta}a \hadam N_\beta(0) -\bar{X}(0)\right\|_\infty + \frac{\hat{a}}{\beta}
  \left (\hat{\lambda} + \hat{\mu} \right ) t \right ) e^{\hat{\mu}
  t}.
$$

Thus, $f_\beta\rightarrow0$ as $\beta\to+\infty$ and
$\sup_{0 \leq s \leq t} \left \| \frac{1}{\beta} a \hadam N_\beta(s) - \bar{X}(s) \right
\|_\infty$ converges in distribution to 0. In particular, for $\varepsilon >0$,

$$
\lim_{\beta \rightarrow +\infty} \Pr \left ( \sup_{0 \leq s \leq t} \left \|
   \frac{1}{\beta} a \hadam N_\beta(s) - \bar{X}(s) \right \|_\infty \geq \varepsilon \right )
= 0.
$$

Hence, $(\beta^{-1} a \hadam N_\beta(t))$ converges to $(\bar{X}(t))$ in probability
uniformly on compact sets.
\ep

The processes $(\beta^{-1} a \hadam N_\beta(t))$ converge in distribution, as $\beta$
tends to infinity, to a process $(\bar{X}(t))$ which is completely
deterministic. This process has a fixed point if and only if there
exists $x\in\R_+^K$ such that
\begin{equation}
  \lambda_k - \mu_k \psi_k(x) = 0\quad \text{for } 1 \leq k \leq K.
  \label{eq:det-stat}
\end{equation}

\begin{proposition}
Consider an acyclic network with $L$ links and $K$ classes.

$(\bar{X}(t))$ admits a unique fixed point if the optimal stability conditions
\eqref{eq:optimal} are satisfied:
$$
\sum_{k:l\in r_k} \rho_k < C_l\ \mathrm{for}\ 1\leq l \leq L.
$$
Conversely, if there exists a link $l_0$ such that
$$
\sum_{k:l_0\in r_k} \rho_k > C_{l_0}
$$
then $(\bar{X}(t))$ does not admit a fixed point.
\end{proposition}
\bp

By definition of acyclic networks, we suppose that the links are numbered in a
such a way that the links on a route of any class are an increasing sequence.
By definition, any fixed point $x=(x_1,\dots,x_K)$ satisfies
$$
\psi_k(x) = \rho_k,\quad \text{for } 1\leq k \leq K.
$$

We first suppose that the optimal stability conditions \eqref{eq:optimal} are satisfied.
First, note that if $x=(x_1,\dots,x_K)$ satisfies
the optimal stability conditions, then for all $k$, $\psi_k(x) = x_k$.
The traffic intensities $(\rho_1,\dots,\rho_K)$ are then a fixed point and this is
the only fixed point which statisfies the optimal stability conditions.

On the contrary, if $x$ is such that there exists $l_0$ with
$$
\sum_{k:l_0\in r_k} x_k \geq C_{l_0},
$$
then, there exists $l_1$ which is saturated i.e.
$$
\sum_{k:l_1\in r_k} \theta_k^{l_1}(x) = C_{l_1}.
$$
Since the network is acyclic, without loss of generality, we can assume that $l_1$ is such that, for all $l>l_1$,
link $l$ is not saturated i.e.
$$
\sum_{k:l\in r_k} \theta_k^l(x) < C_l.
$$
In that case, for all class $k$ such that $l \in r_k$,
we have $\psi_k(x) = \theta_k^{l_1}(x)$ and
\begin{align*}
\sum_{k:l_1\in r_k} \psi_k(x) &= C_{l_1}\\
&> \sum_{k:l_1\in r_k} \rho_k.
\end{align*}
The vector $x$ cannot be a fixed point and the uniqueness of the fixed point is proved.

We now suppose that there exists a link $l_0$ such that
$$
\sum_{k:l_0\in r_k} \rho_k > C_{l_0}
$$
In that case, it is enough to remark that, for all $x=(x_1,\dots,x_K)$,
\begin{align*}
\sum_{k:l_0\in r_k} \psi_k(x) &\leq C_{l_0}\\
&< \sum_{k:l_0\in r_k} \rho_k.
\end{align*}
There is no fixed point.
\ep

This proposition brings the intuition that the stability conditions for any acyclic network
can be arbitrarily close to the optimal ones if the access rates are small enough.
In the remainder  of this section, we prove that this intuition
is true in the case of linear networks and upstream trees.

\subsection{Asymptotic stability of linear networks}
\label{subsection:line-asym}

We return here to the linear networks discussed in Section \ref{section:line}.
Theorems \ref{theo:stability} and \ref{theo:unstability} show that
stability conditions depend on the function $\bar{\phi}_0$ and consequently on stationary
distributions $(\tilde{\pi}^\alpha)$.  In order to characterize the ergodicity
conditions, we have to study the evolution of these distributions when the access
rates decrease to $0$.

Consider the process $\quasi$ introduced in
Section  \ref{subsection:frozen} for $\alpha \in [0,1]$ under the scaling defined in
Section \ref{subsection:scaling-access}. We call $(\quasiscaled{t})$ the
scaled process whose transition rates are given by, for $2\leq k \leq L$,
\begin{align*}
n_k \to n_k+1 &: \beta \lambda_k,\\
n_k \to n_k-1 &: \beta \mu_k \tilde{\psi}_k(n \hadam \beta^{-1}a)
\end{align*}
where the functions $\tilde{\psi}_k$ are defined by Equation \eqref{eq:phitilde}.
We suppose that the ergodicity conditions of this
process are fulfilled, i.e., $\rho_k < 1$ for $2\leq k \leq L$, and denote by
$\tilde{\pi}_\beta^\alpha$ the stationary distribution of
$(\beta^{-1} a_{2:L} \hadam \quasiscaled{t})$.

As in Section  \ref{subsection:scaling-access}, when $\beta$ tends to infinity,
the processes $(\beta^{-1} a_{2:L} \hadam \quasiscaled{t})$ converge in distribution
to a deterministic limit that we call $(\bar{X}^\alpha(t))$.  This
process has a unique fixed point that we call $\gamma(\alpha)$ which is the
solution of the following equations:
$$
\lambda_k-\mu_k\tilde{\psi}_k(\alpha,\vect{x}{2}{L})=0\ \mathrm{for}\ 2\leq k \leq L.
$$
Using the definition of $\tilde{\psi}_k$, $\gamma(\alpha)$ can be built
recursively yielding:
\begin{equation}
  \gamma_k(\alpha) =
  \max \left ( \rho_k, \frac{\rho_k}{1-\rho_k}
    \min \left ( \alpha,   \min_{2\leq j \leq k-1} \left ( 1-\rho_j \right ) \right )  \right )\
  \mathrm{for}\ 2\leq k \leq L.
  \label{eq:fixed-point}
\end{equation}

The next proposition states the convergence of the stationary distributions $\tilde{\pi}^\alpha_\beta$ to
the Dirac measure at $\gamma(\alpha)$ for all $\alpha\in[0,1]$ and the convergence
of $\bar{\phi}_0$.

\begin{proposition}
  Consider a linear network with $L$ links in the quasi stationary regime.
  The set $\{ \tilde{\pi}^\alpha_\beta,\ \alpha \in [0,1],\ \beta \in
  \N \}$ is tight. For any $\alpha\in[0,1]$, the stationary distribution
  $(\tilde{\pi}^\alpha_\beta)$ converges to $\delta_{\gamma(\alpha)}$
  when $\beta\to+\infty$.

  If $\bar{\phi}_{0,\beta}$ is the function such that
  \begin{equation}
    \bar{\phi}_{0,\beta}(\alpha) = \int_{\R^{L-1}} \tilde{\psi}_0(\alpha,\vect{x}{2}{L})
    \tilde{\pi}^\alpha_\beta(\diff\vect{x}{2}{L}),\quad  \text{for  } \alpha \in[0,1],
    \label{eq:phibeta}
  \end{equation}
  then $\bar{\phi}_{0,\beta}$ converges pointwise on $[0,1]$ to
  \begin{equation}
    \bar{\phi}_{0,\infty}:\alpha \mapsto \min \left ( \alpha, \min_{2\leq k \leq L} (1-\rho_k)\right),\quad \text{when } \beta\to+\infty.
  \end{equation}
  Moreover, for all $\alpha\in[0,1]$,
  $$
    \lim_{\beta \to +\infty} \inf_{\alpha\leq u \leq 1} \bar{\phi}_{0,\beta}(u) = \bar{\phi}_{0,\infty }(\alpha).
  $$
\label{prop:line-tight}
\end{proposition}
The tightness and the convergence mentioned in this proposition is the convergence of probability
distributions on $\R_+^{K-2}$ with the usual topology.
\bp
First, we have to prove that the set $\{ \tilde{\pi}^\alpha_\beta,\
\alpha \in [0,1],\ \beta \in \N \}$ is tight.

Note that for $\alpha \in [0,1]$ and for $k\geq2$,
$$
\tilde{\psi}_k(\alpha,\beta\vect{x}{2}{L}) \geq
\frac{n_k a_k}{\beta+ n_k a_k}.
$$
Define the Markov processes $(\hat{N}_{2:L,\beta}(t))$ with transition rates:
\begin{align*}
  n_k\rightarrow n_k+1 &: \lambda_k\beta,\\
  n_k\rightarrow n_k-1 &: \mu_k\beta
  \frac{n_ka_k}{\beta+ n_ka_k}.
\end{align*}
We assume that $\hat{N}_{2:L,\beta}(0)=\tilde{N}^\alpha_{2:L,\beta}(0)$.
We then have $(\beta^{-1} a_{2:L}\hadam \hat{N}_{2:L,\beta}(t)) \geqst (\beta^{-1} a_{2:L} \hadam \tilde{N}^\alpha_{2:L,\beta}(t))$.
In order to prove that the set $\{ \tilde{\pi}^\alpha_\beta,\ \alpha
\in [0,1],\ \beta \in \N \}$ is relatively compact, it is enough to
prove that $\{ \hat{\pi}_\beta,\ \beta \in\N \}$ is tight where
$\hat{\pi}_\beta$ is the invariant measure of
$(\beta^{-1} a_{2:L} \hadam \hat{N}_{2:L,\beta}(t))$. These are birth-death processes whose
components are independent. Their invariant measures thus have a
product form and it is sufficient to study the one-dimensional case:
$$
\hat{\pi}_\beta(n_2\beta a_2) = (1-\rho_2)^{\beta a_2^{-1}+1}
\rho_2^{n_2} \frac{\Gamma\left(\frac{\beta}{a_2}+1 +n_2\right)}{
  \Gamma\left(\frac{\beta}{a_2}+1\right)n_2!}
$$
where $\Gamma$ is the usual gamma function
$$
\Gamma(z) = \int_0^{+\infty}t^{z-1}e^{-t}\diff t,\quad \text{for } z>0.
$$

$\hat{\pi}_\beta$ is a negative binomial distribution and we have
$$
\hat{\pi}_\beta([A,+\infty[) = \frac{\int_{\rho_2}^1
  t^{\beta a_2^{-1}} (1-t)^{\lfloor  \beta A/a_2 \rfloor} \diff t } {
  \int_{0}^1 t^{\beta a_2^{-1}} (1-t)^{\lfloor \beta A/a_2
    \rfloor} \diff t } .
$$
Thus, if $A > \rho_2/(1-\rho_2)$,
\begin{equation}
  \lim_{\beta \rightarrow +\infty} \hat{\pi}_\beta([A,+\infty[) = 0.
  \label{eq:tight}
\end{equation}
Fix $\varepsilon>0$ and $A > \rho_2/(1-\rho_2)$. According to \eqref{eq:tight},
there exists $\beta_0\in\N$ such that for all $\beta\geq \beta_0$,
$\hat{\pi}_\beta([0,A]) \geq 1- \varepsilon$. Additionally, there exists a compact
$\Gamma_1 \subset \R$ such that for $\beta < \beta_0$, $\hat{\pi}_\beta(\Gamma_1) \geq 1- \varepsilon$.
By choosing $\Gamma= [0,A] \cup \Gamma_1$, we have $\hat{\pi}_\beta(\Gamma) \geq 1- \varepsilon$
for all $\beta \in \N$ and the set $\{ \hat{\pi}_\beta, \beta \in\N \}$ is tight on $\R_+$ with the usual topology. Finally,
the set $\{ \tilde{\pi}^\alpha_\beta,\ \alpha \in [0,1],\ \beta \in
\N \}$ is tight  and, thus, relatively compact.

We now prove the convergence of the invariant measures.
The infinitesimal generator $\tilde{\Omega}^\alpha_\beta$ of
$(\beta^{-1} a_{2:L} \hadam \quasiscaled{t})$ is defined by
\begin{align*}
\tilde{\Omega}^\alpha_\beta(f)(\vect{x}{2}{L}) = &\sum_{k=2}^L
 \beta\lambda_k \left(f\left(\vect{x}{2}{L}+\frac{a_k}{\beta} e_k\right)-f\left(\vect{x}{2}{L}\right)\right)\\
 &+ \sum_{k=2}^L \beta\mu_k\tilde{\psi}_k(\alpha,\vect{x}{2}{L}) \left(f\left(\vect{x}{2}{L}-\frac{a_k}{\beta} e_k\right) - f\left(\vect{x}{2}{L}\right)\right).
\end{align*}
The ``infinitesimal generator'' $\bar{\Omega}^\alpha$ of $(\bar{X}^\alpha(t))$
is defined by
$$
\bar{\Omega}^\alpha(f)(\vect{x}{2}{L}) =
\sum_{k=2}^L \lambda_k  \frac{\partial f}{\partial x_k}(\vect{x}{2}{L}) -
\sum_{k=2}^L \mu_k \tilde{\psi}_k(\alpha,\vect{x}{2}{L}) \frac{\partial f}{\partial x_k}(\vect{x}{2}{L}).
$$

Let $f$ be a bounded smooth function  with a bounded derivative.
Since $\{\tilde{\pi}_\beta,\ \beta\in]0,1],\ \alpha \in[0,1]\}$ is tight, for any
$\varepsilon>0$, there exists a compact set $\Gamma \subset \R^{K-2}$ such that
for all $\beta$ in $\N$ and $\alpha$ in $[0,1]$,
\begin{equation}
  \int_{\R^{K-2}_+\setminus \Gamma} f(\vect{x}{2}{L})
  \tilde{\pi}^\alpha_\beta(\diff\vect{x}{2}{L}) \leq  \varepsilon.
  \label{eq:borders}
\end{equation}

Since $f$ and its derivative are bounded and the $\tilde{\psi}_k$ are
$\gamma$-Lipschitz for some $\gamma$, there exists $\eta$ such that for $x$ in
$\Gamma$,
\begin{equation}
  |\tilde{\Omega}^\alpha_\beta(f)(\vect{x}{2}{L}) - \bar{\Omega}^\alpha(f)(\vect{x}{2}{L})|
  \leq  \frac{\eta}{\beta}.
  \label{eq:smooth-f}
\end{equation}
The equilibrium equations give, for all $\alpha$,
\begin{equation}
  \int_{\R_+^{K-2}} \tilde{\Omega}^\alpha_\beta(f)(\vect{x}{2}{L})
  \tilde{\pi}^\alpha_\beta(\diff\vect{x}{2}{L}) = 0.
  \label{eq:gen-equ}
\end{equation}
Since $\{\tilde{\pi}^\alpha_\beta,\ \beta \in \N,\ \alpha \in [0,1] \}$
is relatively compact, we can extract a convergent subsequence
$(\tilde{\pi}^{\alpha_i}_{\beta_i})$ such that $\alpha_i \to \alpha$ and $\beta_i\to+\infty$.
Thanks to relations \eqref{eq:borders}, \eqref{eq:smooth-f}, \eqref{eq:gen-equ}
and the continuity of $f$ and its derivative, we can write that
$$
\lim_{i \rightarrow +\infty} \int_{\R_+^{K-2}} \tilde{\Omega}^{\alpha_i}_{\beta_i}(f)(\vect{x}{2}{L})
\tilde{\pi}^{\alpha_i}_{\beta_i}(\diff\vect{x}{2}{L}) = \int_{\R_+^{K-2}} \bar{\Omega}^\alpha(f)(\vect{x}{2}{L})\bar{\pi}^\alpha(\diff\vect{x}{2}{L}).
$$

Finally, $(\pi^{\alpha_i}_\beta)$ converges to $\bar{\pi}^\alpha$
when $\beta\to+\infty$.

The convergence pointwise of $\bar{\phi}_{0,\beta}$ follows by choosing $\alpha_i=\alpha$.

Finally, for the last part of the proposition, we consider
a sequence $(\alpha_i,\beta_i)$ such that $\alpha_i \to \alpha_\infty$, $\beta_i \to +\infty$
when $i\to\infty$ and, for all $i$, $\bar{\phi}_{0,\beta_i}(u_i)=\inf_{\alpha\leq u \leq 1} \bar{\phi}_{0,\beta_i}(u)$.
We have that
$$
\lim_{i\to \infty} \bar{\phi}_{0,\beta_i}(\alpha_i) = \bar{\phi}_{0,\infty}(\alpha_\infty).
$$
By construction of $\alpha_i$, we also have
$$
\lim_{i\to \infty} \bar{\phi}_{0,\beta_i}(\alpha_i) \leq \inf_{\alpha\leq u \leq 1} \bar{\phi}_{0,\infty}(u).
$$
Since $\bar{\phi}_{0,0}$ is a non-decreasing function, its implies that
$$
\lim_{i\to \infty} \bar{\phi}_{0,\beta_i}(\alpha_i) = \bar{\phi}_{0,\infty}(\alpha).
$$
\ep

We finally deduce the following theorem which is a consequence of
Theorem \ref{theo:stability} and Proposition \ref{prop:line-tight}.

\begin{theorem}
  In a linear network with $L$ links, for any traffic intensities
  $(\rho_0,\dots,\rho_L)$ satisfying the optimal stability conditions \eqref{eq:line-optimal},
  if the access rates $(a_0,\dots,a_L)$ are small enough,
  the resulting stochastic process is ergodic.
\end{theorem}
\bp
We consider the scaled version of the process introduced in Subsection \ref{subsection:scaling-access}
$(N_\beta(t))$. According to Theorem \ref{theo:stability}, a sufficient condition
of stability is
\begin{align*}
\rho_k &< 1,\quad \text{for } 0\leq k\leq L,\\
\rho_0 &< \inf_{1-\rho_1\leq \alpha \leq 1} \bar{\phi}_{0,\beta}(\alpha).
\end{align*}
Since traffic intensities $(\rho_0,\dots,\rho_L)$ satisfy optimal stability conditions \eqref{eq:line-optimal},
there exists $\varepsilon>0$ such that
\begin{align*}
\rho_k &< 1,\quad \text{for } 0\leq k\leq L,\\
\rho_0 &< \min_{1\leq k \leq L} (1-\rho_k) - \varepsilon.
\end{align*}
According to Proposition \ref{prop:line-tight}, there exists $\beta_0$ such that
for all $\beta \geq \beta_0$, we have
$$
\inf_{1-\rho_1\leq \alpha \leq 1} \bar{\phi}_{0,\beta}(\alpha) \geq \min_{1\leq k \leq L} (1-\rho_k) - \varepsilon/2.
$$
This implies that, for all $\beta\geq \beta_0$, the process $(N_{\beta}(t))$ is
ergodic.

Finally, to conclude, it is enough to remark that the process
with transition rates
\begin{align*}
n_k \to n_k+1 &:\quad \lambda_k,\\
n_k \to n_k-1 &:\quad \mu_k\psi_k(n\hadam \beta^{-1} a)
\end{align*}
admits the same stationary measures as $(N_{\beta}(t))$ and is then ergodic
if and only if $(N_{\beta}(t))$ is ergodic.
\ep

\begin{figure}[!h]
  \centering
  \begin{minipage}[t]{.4\linewidth}
    \centering
    \includegraphics[width=.95\linewidth]{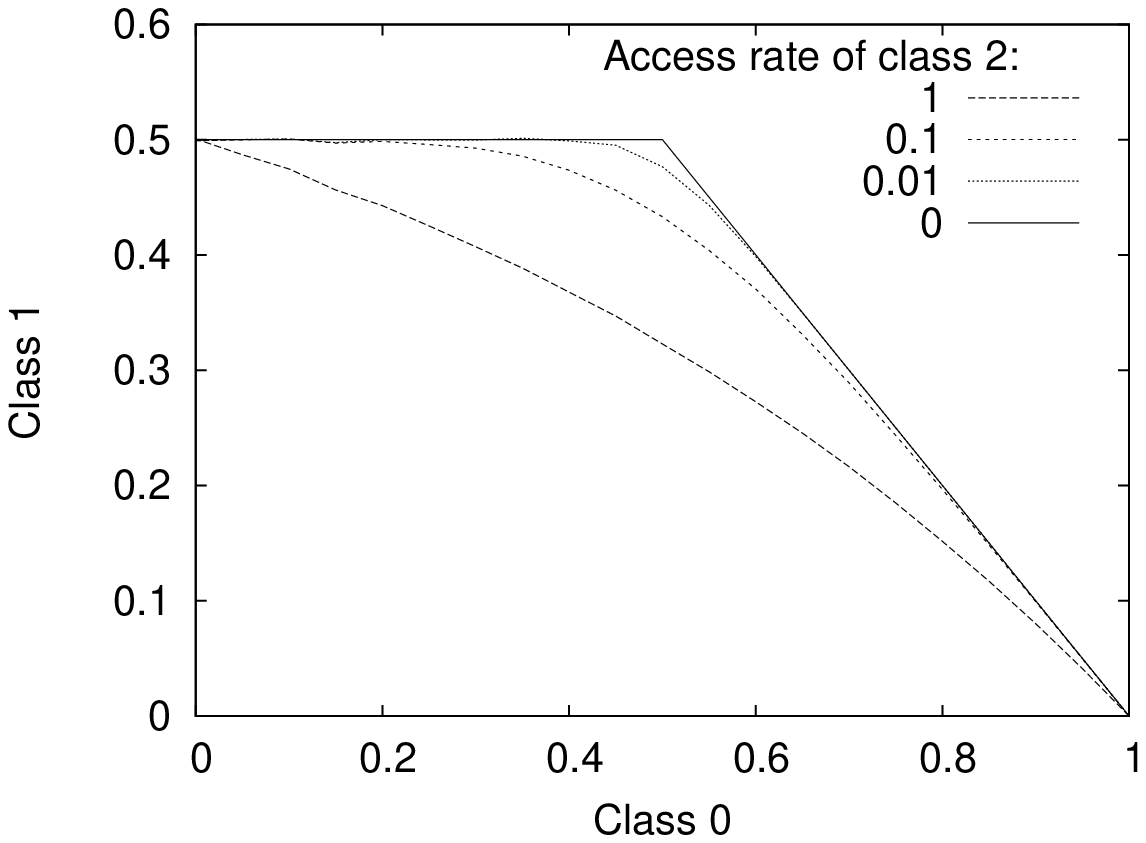}
    \caption{Stability region for classes $0$ and $1$ when $L=2$ and
      $\rho_2=0.5$}
    \label{fig:stability-line-2}
  \end{minipage}
  \hspace{.05\linewidth}
  \begin{minipage}[t]{.4\linewidth}
    \centering
    \includegraphics[width=.95\linewidth]{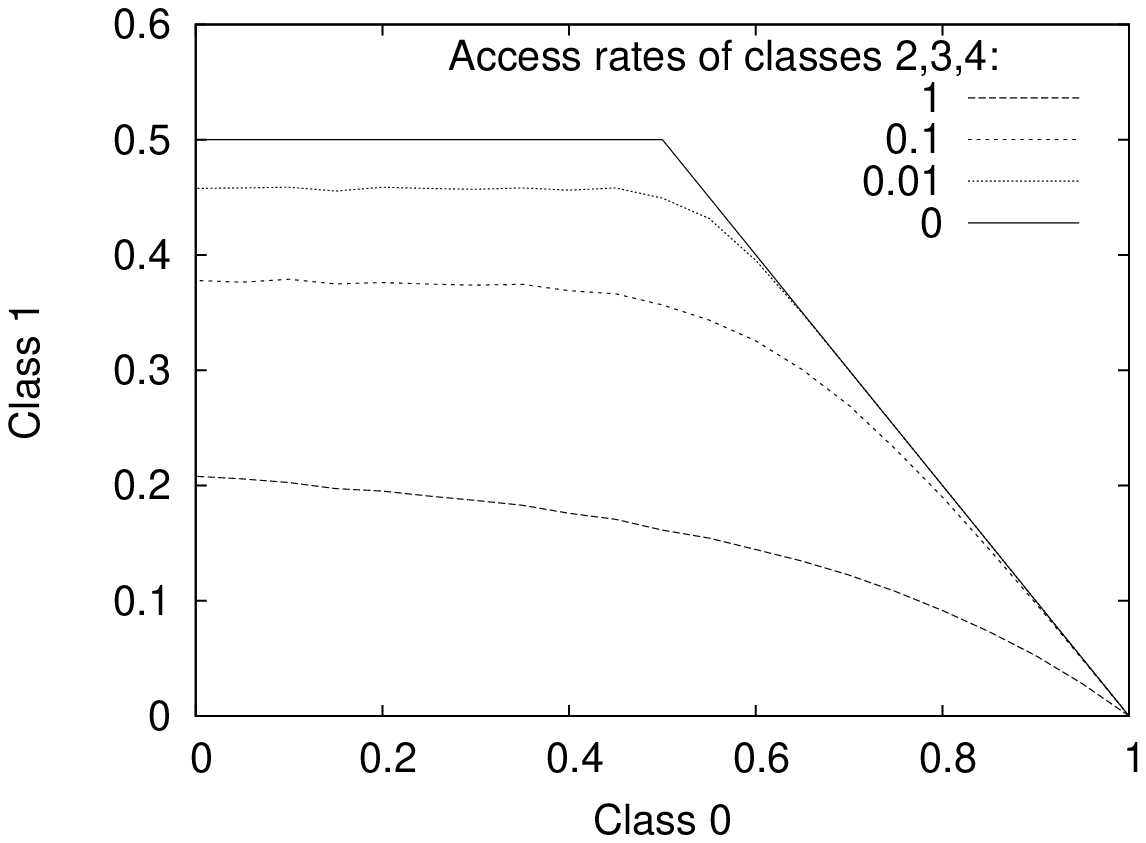}
    \caption{Stability region for classes $0$ and $1$ when $L=4$ and
      $\rho_2,\rho_3,\rho_4=0.5$}
    \label{fig:stability-line-4}
  \end{minipage}
\end{figure}

This result means that the Tail Dropping policy is asymptotically
optimal in linear networks.  Figures \ref{fig:stability-line-2} and
\ref{fig:stability-line-4} represent the stability region for classes
$0$ and $1$, obtained by simulation, for two particular networks and illustrate this result.

\subsection{Asymptotic stability of upstream trees}
\label{subsection:tree-asym}

Upstream trees are a specific class of networks because they are monotonic and this
will allow us to use stochastic domination to establish the asymptotic
stability of upstream trees. Thanks to the monotonicity, the worst case, in the sense of stochastic domination,
for a given class is when there is an infinity of flows in all other classes.
When the optimal stability conditions are satisfied, we can prove that there
exists a class $k_0$ which can get enough bandwidth to be stable even if there is
an infinity of flows in all other classes. Thanks to the same scaling that we use
in subsections \ref{subsection:scaling-access} and  \ref{subsection:line-asym},
we are able to quantify the bandwidth used by the class $k_0$ when the access
rates are decreasing to $0$ and to prove that there exists a class $k_1$
which is stable when the access rates are small enough even if there
is an infinity of flows in all classes except $k_0$ and $k_1$. We conclude by
using a recursion.

In upstream trees,
when two different classes go through the same link, they
share the same path until they exit the network. The last link
on the path of all classes is the same and it is called the root.
Formally, upstream trees are defined as follows:
\begin{description}
\item[(i)] a common root link, say 1, so that $r_k(d_k)=1$ for all
  classes $k$;

\item[(ii)] for any two classes $j$ and $k$, there exists $m$ such that they only have
their last $m$ links in common; i.e, for $i$ in
$\{0,\dots,m-1\}$, $r_j(d_j-i) = r_k(d_k-i)$ and for $i_k$ in $\{0,\dots,d_k-m\}$ and
$i_j$ in $\{0,\dots,d_j-m\}$, $r_j(i_j) \neq r_k(i_k)$.
\end{description}

\begin{figure}[h]
  \centering
  \begin{minipage}{0.3\linewidth}
    \centering
    \includegraphics[width=.8\linewidth]{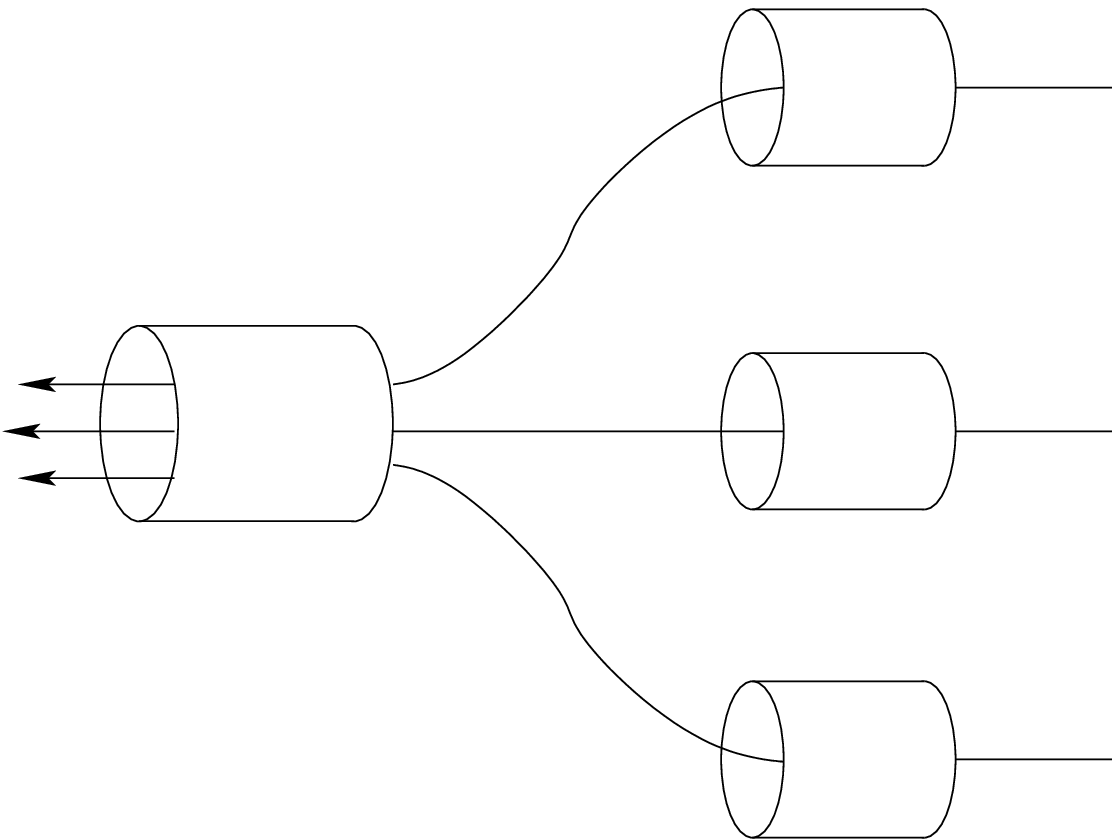}
  \end{minipage}
  \begin{minipage}{.3\linewidth}
    \centering
    \includegraphics[width=.8\linewidth]{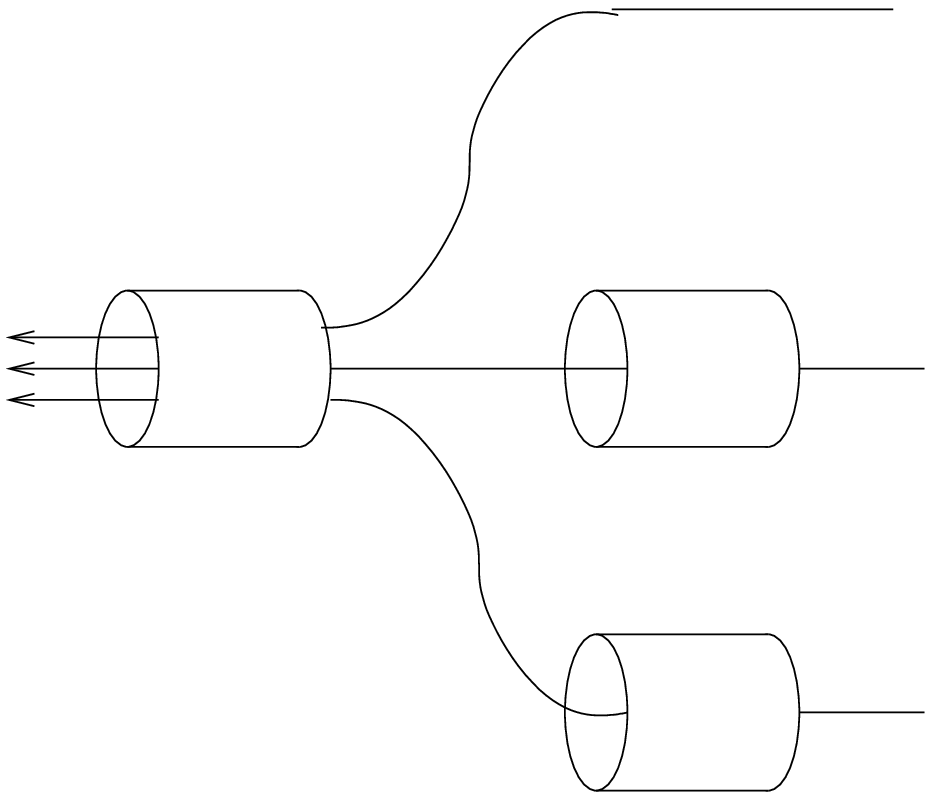}
  \end{minipage}
  \caption{Examples of upstream trees}
  \label{fig:uptrees}
\end{figure}

Two examples of upstream trees are given in Figure
\ref{fig:uptrees}.  In \cite{bonald-09}, two very specific
examples were shown not to be ergodic under optimal stability conditions except
when the access rates tend to $0$. We generalize this result.

For any link $l$, let $S_l$ denote the set of children of $l$, the links
just \emph{before} $l$ on the path of some class:
$$
S_l = \bigl\{ j\in\{1,\dots,L\},\ \exists k \in\{1,\dots,K\},\ \exists i,\ r_k(i)=j,
r_k(i+1)=l \bigr\}.
$$
In particular, $S_1$ is the set of links which are the links just before the
root link on the path of all classes with a route of length longer than $1$.

We suppose that there are no two classes with exactly the same path. We also suppose
without loss of generality that all the
links of the network can be saturated. For any link $l$, this implies that there
is a class $k_l$ directly entering in the network through that link, i.e.,
$r_{k_l}(1)=l$ or $\sum_{j \in S_l} C_j > C_l$. If a link $l_0$ cannot be
saturated then the bandwidth allocation is the same with or without $l_0$ and we can remove it.

Upstream trees are monotonic in the sense that if $n\in\N^K$ and
$m\in\N^K$ and $m_k=n_k$ for a given $k$ and $m_j\geq n_j$ for $j\neq
k$ then $\phi_k(m) \leq \phi_k(n)$. Thanks to this monotonicity property, the
processes we study fit into the framework of \cite{borst-08}. We suppose that
the traffic intensities satisfy the optimal stability conditions
\eqref{eq:optimal}.

Because of the monotonicity of this network, we know that the worst case for class
$k_0$ is when the numbers of flows in all other classes increase to infinity.
We define the following allocation:
\begin{equation}
\forall x_{k_0}\in \R_+,\ \satpsi{1}{k_0}(x_{k_0}) = \lim_{\eta \rightarrow \infty}
\inf_{\substack{y:y_{k_0}=x_{k_0},\\ k\neq k_0, y_k>\eta}} \psi_{k_0}(y)
\label{eq:sat-psi-1}
\end{equation}
and $\satphi{1}{k_0}$ is such that $\satphi{1}{k_0}(n_{k_0})=
\satpsi{1}{k_0}(a_{k_0}n_{k_0})$ for all $n_{k_0} \in \N$.
We then define the Markov process $(\satN{1}{k_0}(t))$ such that
$\satN{1}{k_0}(0)=N_{k_0}(0)$ with the transition rates
\begin{align*}
n_{k_0} \rightarrow n_{k_0}+1&: \lambda_{k_0},\\
n_{k_0} \rightarrow n_{k_0}-1&: \mu_{k_0}\satphi{1}{k_0}(n_{k_0}).
\end{align*}
Because of the monotonicity of the considered allocations, we have
that the process $(\satN{1}{k_0}(t))$ stochastically dominates $(N_{k_0}(t))$ (see \cite{borst-08}).

In the sense of the stochastic domination, this process represents the worst
case for class $k_0$. The next lemma shows that there is a class which is stable
even in the worst case.

\begin{lemma}
Consider an upstream tree and suppose the traffic intensities of this network
respect the optimal stability conditions \eqref{eq:optimal}.
There exists a class $k_0$ such that the Markov process $(\satN{1}{k_0}(t))$
is ergodic.
\label{lemma:ergo-k0}
\end{lemma}
\bp
Because of the monotonicity of the network all that we have to prove is that
there exists a class $k_0$ such that
  \begin{equation}
    \lim_{\eta \rightarrow +\infty}
    \inf_{n_1\geq\eta,\dots,n_K \geq\eta} \phi_{k_0}(n) > \rho_{k_0}.
    \label{eq:up-k-stable}
  \end{equation}
If $k_0$ satisfies \eqref{eq:up-k-stable} then there exists $\varepsilon>0$
and $\eta_0$ such that, for all $n_{k_0}\geq \eta_0$,
$$
\lambda_{k_0}-\mu_{k_0}\satphi{1}{k_0}(n_{k_0}) \leq -\varepsilon.
$$
This implies the ergodicity of $(\satN{1}{k_0}(t))$. We proceed by recursion on
the depth of the tree.

If the depth is $1$ then there is a single class and a single link and
condition \eqref{eq:up-k-stable} is obviously true.

We suppose now that the depth of the tree is $d\geq2$. If there is a class $k_0$
such that $r_{k_0}(1)=1$ then the root link
is the entry point of $k_0$ which is the only one entering the network through $1$. The worst case
for class $k_0$ is when all links in $S_1$ are saturated; we have for all
$n\in \N^K$,
$$
\phi_{k_0}(n) \geq \min\left(a_{k_0}n_{k_0},
C_1\frac{a_{k_0}n_{k_0}}{ \sum_{l\in S_1} C_l +a_{k_0}n_{k_0}}\right).
$$
Because the optimal conditions are satisfied, we have $\rho_{k_0}<C_1$ and there exist
$\delta>0$ and $\eta_0$ such that for all $n\in\N^K$ with $n_{k_0}\geq\eta_0$, we
have
$$
\phi_{k_0}(n) \geq \rho_{k_0}+\delta
$$
which is enough to conclude.

We now suppose that there is no class directly connected to the root,
i.e. there is no class $k$ with a path of length $1$. This implies
that $S_1$ is not empty. Since the root link can be saturated, we have
$\sum_{l\in S_1} C_l > C_1$.
Due to optimal stability conditions, we have
$$
\sum_{k=1}^K \rho_k < C_1
$$
and we can deduce that there is a link $l_0$ in $S_1$ such that
$$
\sum_{k:r_k(d_k-1)=l_0} \rho_k < \frac{C_{l_0}}{\sum_{j\in S_1}C_j}C_1.
$$
Since the link $l_0$ can be saturated, there is a constant $\eta_0$
such that, if for any class $k$ going through $l_0$ (i.e. satisfying $r_k(d_k-1)=l_0$),
we have $n_k>\eta_0$, we have
$$
\sum_{k:r_k(d_k-1)=l_0} \theta_k^{d_k-1}(n\hadam a) = C_{l_0}.
$$
Using the monotonicity of upstream trees, we can deduce that the worst
case for classes going through $l_0$ is when all the other links in $S_1$
are saturated. We deduce immediately that when $n_k>\eta_0$ for all $k$
going through $l_0$, we have
$$
\phi_k(n) \geq
\frac{C_1}{\sum_{j\in S_1}C_j}\theta_k^{d_k-1}(n\hadam a).
$$
It remains to consider the case where all the links in $S_1$
are saturated corresponding to the case of equality in the previous equation.
With these conditions, the subtree whose root link is $l_0$ is equivalent
to an upstream tree where the capacity of link $l_0$ is replaced by
$C_1C_{l_0}/(\sum_{j\in S_1} C_j)$.

If there is a class $k_0$ such that $d_{k_0}=l_0$ then, as previously, we have
$$
\lim_{\eta \rightarrow \infty}\inf_{n:n_k\geq\eta} \phi_{k_0}(n) \geq
\frac{C_1C_{l_0}}{\sum_{j\in S_1} C_j} > \rho_{k_0}
$$
and $k_0$ satisfies \eqref{eq:up-k-stable}.

If there is no class $k_0$ such that $d_{k_0}=l_0$ then $S_{l_0}$ is not empty
and, as previously, we can find $l_1$ such that
$$
\sum_{k:r_k(d_k-2)=l_1} \rho_k < C_{l_1}\frac{C_{l_0}}{\sum_{j\in S_{l_0}} C_j} \frac{C_1}{\sum_{j\in S_1}C_j}.
$$

Using a recursion, we finally find $k_0$ satisfying \eqref{eq:up-k-stable}.
In addition, the links on the path of $k_0$ satisfy
\begin{equation}
\forall i\in\{1,\dots,d_{k_0}-1\},\ \sum_{k:r_{k_0}(i) \in r_k} \rho_k <
C_{r_{k_0}(i)} \prod_{j=i+1}^{d_{k_0}} \frac{C_{r_{k_0}(j)}}{\sum_{l\in S_{r_{k_0}(j)}} C_l}.
\label{eq:succ-links}
\end{equation}

\ep

We now perform the scaling of Section \ref{subsection:scaling-access} on the
process $(\satN{1}{k_0}(t))$ and call the scaled process
$(\satscaledN{1}{k_0}(t))$. We deduce from the above lemma that
the scaled process is ergodic for all $\beta\in\N$. For each $\beta$, we call
$\satpi{1}{k_0}$ the stationary distribution of $(\beta^{-1} a_{k_0}\satscaledN{1}{k_0}(t))$.
As in Section \ref{subsection:scaling-access}, the process
$(\beta^{-1} a_{k_0}\satscaledN{1}{k_0}(t))$ converges in probability uniformly on compact sets
to a deterministic process
$(\satX{1}{k_0}(t))$. The next lemma shows that there is also convergence
for the stationary distribution. The proof of this lemma is similar to that
of Proposition \ref{prop:line-tight} and is omitted. The tightness and convergence
in distribution mentioned in this lemma are on the space $\R_+$ with its usual topology.

\begin{lemma}
Consider an upstream tree and suppose that the optimal stability conditions \eqref{eq:optimal}
are satisfied by the traffic intensities. Consider $k_0$ such that $(\satN{1}{k_0}(t))$
is ergodic. The set $\{ \satpi{1}{k_0}, \beta\in\N \}$
is tight and $\satpi{1}{k_0}$ converges to $\delta_{\alpha^1_{k_0}}$ where
$\alpha^1_{k_0}$ is the unique solution of the following equation
\begin{equation}
\lambda_{k_0} - \mu_{k_0} \satpsi{1}{k_0}(\alpha^1_{k_0})=0.
\label{eq:sat-fix-point-1}
\end{equation}
The function $\satpsi{1}{k_0}$ is defined by Equation \eqref{eq:sat-psi-1}.
\label{lemma:conv-k0}
\end{lemma}

Now, we proceed by recursion. First, we exhibit a class $k_1$ which is
also stable when the access rates of classes $k_0$ and $k_1$ are small enough whatever
the state of other classes. For any classes $k_0$ and $k_1$, define the
following allocation
\begin{equation}
\forall (x_{k_0},x_{k_1})\in\R_+^2,\
\satpsi{2}{k_0}(x_{k_0},x_{k_1}) = \lim_{\eta \rightarrow \infty}
\inf_{\substack{y:y_{k_0}=x_{k_0}\\y_{k_1}=x_{k_1}\\k\neq k_i,y_k\geq\eta}}
\psi_{k_0}(x)
\label{eq:sat-psi-2}
\end{equation}
and $\satpsi{2}{k_1}$ is defined similarly. Now define $\satphi{2}{k_0}$ and
$\satphi{2}{k_1}$ such that $\satphi{2}{k_i}(n_{k_0},n_{k_1}) =
\satpsi{2}{k_i}(a_{k_0}n_{k_0},n_{k_1}n_{k_1})$ for all $(n_{k_0},n_{k_1})$
in $\N^2$. Let $(\satN{2}{k_0}(t),\satN{2}{k_1}(t))$ be the Markov process
with the transition rates for $i=0,1$:
\begin{align*}
n_{k_i} &\rightarrow n_{k_i}+1 : \lambda_{k_i},\\
n_{k_i} &\rightarrow n_{k_i}-1 : \mu_{k_i} \satphi{2}{k_i}(n_{k_0},n_{k_1}),
\end{align*}
and such that $(\satN{2}{k_0}(0),\satN{2}{k_1}(0))=(N_{k_0}(0),N_{k_1}(0))$.
Using the monotonicity of the considered allocations, we know that
process $(\satN{2}{k_0}(t),\satN{2}{k_1}(t))$ stochastically dominates
$(N_{k_0}(t),N_{k_1}(t))$. Again, let $(\satscaledN{2}{k_0}(t),\satscaledN{2}{k_1}(t))$
denote the scaled version of $(\satN{2}{k_0}(t),\satN{2}{k_1}(t))$.

\begin{lemma}
Consider an upstream tree and assume that the optimal stability conditions \eqref{eq:optimal}
are satisfied by the traffic intensities of this tree. Consider $k_0$ as defined
in Lemmas \ref{lemma:ergo-k0} and \ref{lemma:conv-k0}.

There exist $k_1$ and $\beta_1\in\N$ such that, for all $\beta\geq \beta_1$, the
Markov process $(\satscaledN{2}{k_0}(t),\satscaledN{2}{k_1}(t))$ is ergodic.
When $\beta\geq\beta_1$, the Markov process
$$(\beta^{-1} a_{k_0}\satscaledN{2}{k_0}(t), \beta^{-1} a_{k_1}\satscaledN{2}{k_1}(t))$$
has a unique stationary distribution that we denote by $\satpi{2}{k_0,k_1}$.
The set
$$\left\{ \satpi{2}{k_0,k_1},\ \beta \in\N \right\}$$
is tight and the stationary distribution $\satpi{2}{K_0,k_1}$ converges to
$\delta_{(\alpha^2_{k_0},\alpha^2_{k_1})}$ where $(\alpha^2_{k_0},\alpha^2_{k_1})$
is the unique solution of the equations
\begin{equation}
\left \{
\begin{aligned}
\lambda_{k_0} - \mu_{k_0} \satpsi{2}{k_0}\left(\alpha^2_{k_0},\alpha^2_{k_1}\right) &=0,\\
\lambda_{k_1} - \mu_{k_1} \satpsi{2}{k_1}\left(\alpha^2_{k_0},\alpha^2_{k_1}\right) &=0.
\end{aligned}
\right.
\label{eq:sat-fix-point-2}
\end{equation}
The functions $\satpsi{2}{k_0}$ and $\satpsi{2}{k_1}$ are defined by Equation
\eqref{eq:sat-psi-2}.
\label{lemma:ergo-k1}
\end{lemma}
The tightness and convergence
in distribution mentioned in this lemma are on the space $\R^2_+$ with its usual topology.

\bp
First, we prove that there exists $k_1$ such that
\begin{equation}
\lim_{x_{k_1}\rightarrow\infty} \satpsi{2}{k_1}(\alpha_{k_0}^1,x_{k_1}) > \rho_{k_1}
\label{eq:lim-k1}
\end{equation}
where $\alpha_{k_0}^1$ is the unique solution of \eqref{eq:sat-fix-point-1}.

For all $i$ in $\{1,\dots,d_{k_0}\}$, We define
$$
\breve{\tau}^i_{k_0}\left(x_{k_0}\right) =
\lim_{\eta\rightarrow \infty}
\inf_{\substack{y:y_{k_0}=x_{k_0}\\k\neq k_0,y_k\geq\eta}}
\theta^i_{k_0}(y).
$$
Using Lemma \ref{lemma:conv-k0}, we can prove
$$
\sigma_i = \breve{\tau}^i_{k_0}\left(\alpha_{k_0}^1\right) =
\lim_{\beta\rightarrow +\infty} \E_{\satpi{1}{k_0}}\left(\breve{\tau}^i_{k_0}\left( \beta^{-1} a_{k_0}\satscaledN{1}{k_0} \right)  \right)
$$
The quantity $\sigma_i$ is the bandwidth used asymptotically by class $k_0$
in link $r_{k_0}(i)$ when all the other classes are saturated and when
the access rate $a_{k_0}$ decreases to $0$. We want to prove that when $a_{k_0}$
is small enough, class $k_0$ leaves enough bandwidth for other classes.
First, we have the following result which states that the optimal stability
conditions are met by all classes except $k_0$ when we remove the bandwidth
asymptotically used by class $k_0$ in the saturated case:
$$
\forall i\in\{1,\dots,d_{k_0}\},\ \sum_{k:r_{k_0}(i)\in r_k} \rho_k < C_{r_{k_0}(i)} - \sigma_i.
$$
Indeed, for link 1, note that $\sigma_{d_{k_0}}=\rho_{k_0}$ because $\alpha_{k_0}^1$
is the solution of Equation \eqref{eq:sat-fix-point-1} and the previous inequality is obviously true for
link 1. If the route of $k_0$ is of length 1 the lemma is proved.

Suppose that the route of $k_0$ is at least of length 2. By construction of $k_0$,
there is no class entering the network through $r_{k_0}(i)$ for $i$ in
$\{2,\dots,d_{k_0}\}$. Since any link can be saturated, this implies
$$
\forall i\in \{2,\dots,d_{k_0}\},\  C_{r_{k_0}(i)} < \sum_{l\in S_{r_{k_0}(i)}} C_l.
$$
Moreover, we can prove that
$$
\forall i\in\{1,\dots,d_{k_0}-1\},\  \sigma_i = \prod_{j=i+1}^{d_{k_0}}
\frac{\sum_{l\in S_{r_{k_0}(j)}} C_l}{C_{r_{k_0}}(j)} \rho_{k_0}.
$$
Using these last two equations and Equation \eqref{eq:succ-links}, we deduce
$$
\forall i\in\{1,\dots,d_{k_0}-1\},\ \sum_{k:r_{k_0}(i)\in r_k} \rho_k < C_{r_{k_0}(i)} - \sigma_i.
$$

Considering, for all $k\neq k_0$, the allocation
$$
(x_0,\dots,x_{k_0-1},x_{k_0+1},\dots,x_K) \mapsto
\psi_k\left(x_0,\dots,\alpha_{k_0}^1,\dots,x_K\right)
$$
and using the same procedure used to find $k_0$, we can find $k_1$ such
that
$$
\lim_{\eta \rightarrow \infty} \inf_{\substack{k\neq k_0 x_k \geq \eta\\x_{k_0}=\alpha_{k_0}^1}}
\psi_{k_1}(x)> \rho_{k_1}
$$
which is equivalent to Equation \eqref{eq:lim-k1}.

We are now able to prove that there exists $\beta_1\in\N$ such that
for all $\beta\geq \beta_1$, the Markov process $(\satscaledN{2}{k_0}(t),\satscaledN{2}{k_1}(t))$
is ergodic. By construction of $(\satscaledN{1}{k_0}(t))$ and
$(\satscaledN{2}{k_0}(t),\satscaledN{2}{k_1}(t))$, according to Theorem 2 of \cite{borst-08},
we just have to prove that there exists $\beta_1\in\N$ such that, for all $\beta\geq\beta_1$,
\begin{equation}
\rho_{k_1} < \E_{\satpi{1}{k_0}} \left (
\liminf_{x_{k_1} \rightarrow \infty} \satpsi{2}{k_1}\left( \beta^{-1}a_{k_0} \satscaledN{1}{k_0}(0),x_{k_1}\right)\right).
\label{eq:borst-2}
\end{equation}
We deduce from \eqref{eq:lim-k1} that there exists $b>\alpha_{k_0}^1$,
$\varsigma>0$ and $\eta_0$ such that
\begin{equation}
\forall x_{k_0}\leq b,\ \forall x_{k_1}\geq\eta_0,\ \satpsi{2}{k_1}\left(x_{k_0},x_{k_1}\right)\geq \rho_{k_1}+\varsigma.
\label{eq:blah}
\end{equation}
We know that $\satpi{1}{k_0}$ converges to $\delta_{\alpha_{k_0}^1}$ when
$\beta\to+\infty$. Thus, there exists $\beta_1\in\N$, such that, for all
$\beta\geq\beta_1$,
\begin{equation}
\satpi{1}{k_0}\left( [b, \infty[\right) \leq \frac{\varsigma}{2(\rho_{k_1}+\varsigma)}.
\label{eq:pi2-tight-1}
\end{equation}
Finally, we have that, for all $\beta\geq\beta_1$, Equation \eqref{eq:borst-2} is
satisfied and the process $(\satscaledN{2}{k_0}(t),\satscaledN{2}{k_1}(t))$
is ergodic.

We now prove tightness. By construction, we have
$(\satscaledN{1}{k_0}(t))\geqst(\satscaledN{2}{k_0}(t))$ and we deduce immediately
that
\begin{equation}
\forall  c \in \R,\ \satpi{2}{k_0,k_1}([c,\infty[\times \R) \leq
\satpi{1}{k_0}([c,\infty[).
\label{eq:pi2-tight-2}
\end{equation}
Moreover, according to Equations \eqref{eq:blah} and \eqref{eq:pi2-tight-1},
for all $\beta\geq\beta_1$, there exists a process $(Y_\beta(t))$ such that
$(Y_\beta(t))\geqst(\satscaledN{2}{k_1}(t))$ with the following transition rates
\begin{align*}
n_{k_1} &\rightarrow n_{k_1}+1: \beta\lambda_{k_1},\\
n_{k_1} &\rightarrow n_{k_1}+1: \beta\left(\lambda_{k_1} + \frac{\varsigma}{3}\mu_{k_1}\right)
\ind_{\{Y_\beta(t)\geq \eta_0 \beta\}}.
\end{align*}
$(Y_\beta(t))$ has the transition rates of a $M/M/1$ queue with $\lfloor \eta_0 \beta \rfloor$
permanent clients.
Let $\pi_{Y,\beta}$ denote the stationary distribution of $(\beta^{-1} Y_\beta(t))$ and we
obtain immediately that
\begin{align*}
\forall \beta\geq\beta_1,\ \forall \eta\geq\eta_0,\ \pi_{Y,\beta}([\eta,\infty[) &\leq
\left( \frac{\lambda_{k_1}}{\lambda_{k_1} + \frac{\varsigma}{3}\mu_{k_1}}\right)^{\left\lfloor(\eta-\eta_0)\beta\right\rfloor}\\
&\leq \left( \frac{\lambda_{k_1}}{\lambda_{k_1} + \frac{\varsigma}{3}\mu_{k_1}}\right)^{\left\lfloor(\eta-\eta_0)\beta_1\right\rfloor}
\end{align*}
We deduce from the previous equation that
\begin{equation}
\forall \beta\geq\beta_1,\ \forall \eta\geq\eta_0,\ \satpi{2}{k_0,k_1}(\R \times [c,\infty])
\leq \left( \frac{\lambda_{k_1}}{\lambda_{k_1} + \frac{\varsigma}{3}\mu_{k_1}}\right)^{\left\lfloor(\eta-\eta_0)\beta_1\right\rfloor}
\label{eq:pi2-tight-3}
\end{equation}

The tightness of $\{\satpi{2}{k_0,k_1},\ \beta\geq\beta_1 \}$ follows from
Equations \eqref{eq:pi2-tight-2} and \eqref{eq:pi2-tight-3} and the tightness
of $\{\satpi{1}{k_0},\ \beta\in\N \}$.

We can remark that existence and uniqueness of the solution
of Equation \eqref{eq:sat-fix-point-2} come from the monotonicity of
$(\satpsi{2}{k_0},\satpsi{2}{k_1})$ and from Equations \eqref{eq:up-k-stable} and
\eqref{eq:lim-k1}.

The proof of convergence is similar to that used in the proof of Proposition
\ref{prop:line-tight}.
\ep

Using a recursion we can easily  generalize the previous lemma. The principle
is still the same: little by little, we decrease the access rates and the classes
become stable one by one. In particular, there exists $\beta_{K-1}\in\N$ such that
for all $\beta\geq \beta_{K-1}$, the process $(N_\beta(t))$ is ergodic
where $(N_\beta(t))$ is the scaled version of $(N(t))$.
We can therefore state the following.

\begin{theorem}
  In an upstream tree, for any traffic intensities
  satisfying the optimal stability conditions \eqref{eq:optimal},
  there exist positive access rates  small enough
  such that the resulting stochastic process is ergodic.
\end{theorem}

%% file: proof1.tex
We first give a sketch of the proof. The technical details are proved in the
lemmas which follow.

We consider the process $(\bar{N}(m,t))=((\bar{N}_0(m,t),\bar{N}_1(m,t),\dots,\bar{N}_k(m,t)))$
such that, for all $0\leq k \leq L$,
$$
\bar{N}_k(m,0) = \frac{m_k}{\|m\|}\quad \text{and}\quad
\bar{N}_k(m,t) = \frac{1}{\|m\|} N_k(\|m\|t),\quad \text{for } t \geq 0.
$$
We consider a sequence $(m^i,\ i\in\N)$ such that
\begin{equation}
\begin{aligned}
&\lim_{i\to+\infty} \|m^i\| = +\infty,\quad
&&\lim_{i\to+\infty} \frac{m^i_0}{\|m^i\|} = \alpha,\\
&\lim_{i\to+\infty} \frac{m^i_1}{\|m^i\|} = 1-\alpha,\quad
&&\lim_{i\to+\infty} \frac{m^i_k}{\|m^i\|} = 0,\quad \text{for } 2\leq k \leq L.
\end{aligned}
\label{eq:mi}
\end{equation}
 We want to prove that the sequence of processes $(\bar{N}(m^i,t))$
converge and to characterize the limit. For that purpose, we write them as
the sum of a martingale term and a continuous term and we define
the process $(\bar{M}(m,t))$ such that for $t\geq0$ and
$0\leq k \leq L$, we have
\begin{equation}
\begin{split}
  \bar{M}_k(m,t) = &\bar{N}_k(m,t) - \frac{m_k}{\|m\|} - \lambda_k t\\ &+ \mu_k
  \int_0^t \phi_k\left(\bar{N}_0(m,s),\bar{N}_1(m,s),N_2(\|m\|s),\dots,N_L(\|m\|s)\right)\diff s.
\end{split}
  \label{eq:mbar}
\end{equation}

Because the access rates $a_0$ and $a_1$ are equal to 1, we have
$$
\forall n\in\N^K,\ \gamma \in \R,
\phi_k\left(\gamma n_0,\gamma n_1,n_2,\dots,n_L\right) =
\phi_k\left(n_0,n_1,n_2,\dots,n_L\right)
$$
If $a_0$ and $a_1$ are not equal to $1$, the previous expression is true outside
a compact set and it is possible to prove all the results in this section
without this assumption because, at a fluid level, the first link is always saturated.

In Lemma \ref{lemma:relat-compact}, we prove that the set
$$
\mathcal{C} = \left \{ (\bar{N}(m^i,t)),\ i\in\N \right \}
$$
is relatively compact and its limiting points are continuous. We can now extract a converging subsequence
and we suppose that $(m^i)$ is such that $(\bar{N}(m^i,t))$ converges in distribution
to a limit that we denote $(\bar{Z}(t))=(\bar{Z}_0(t),\bar{Z}_1(t),\dots,\bar{Z}_L(t))$
and which is continuous.
In lemma \ref{lemma:2-K}, we characterize $(\bar{Z}_2(t),\dots,\bar{Z}_L(t))$ by proving
that, for $2\leq k \leq L$ and $t\geq 0$, $\bar{Z}_k(t)=0$.

Proposition \ref{prop:averaging} is the key result in order to understand the time
scale separation between classes 0 and 1 and classes $2,\dots,L$. We deduce
from this lemma that
$$
\left(\int_0^t \phi_k \left(\bar{N}_0(m^i,s),\bar{N}_1(m^i,s),N_2(\|m^i\|s),\dots,N_L(\|m^i\|s)\right)\diff s,\
k=0,1\right)
$$
converges in distribution to
$$
\left(\int_0^t \bar{\phi}_0\left(\frac{\bar{Z}_0(s)}{\bar{Z}_0(s)+\bar{Z}_1(s)}\right)\diff s,\
\int_0^t \frac{\bar{Z}_1(s)}{\bar{Z}_0(s)+\bar{Z}_1(s)}\diff s\right)
$$
where $\bar{\phi}_0$ characterizes the average throughput of class $0$ in the
quasi-stationary case as defined by Equation \eqref{eq:phi0}.

We can deduce from that and the fact that $(\bar{M}(m^i,t))$ converges to 0 in
distribution that $(\bar{N}_0(m^i,t),\bar{N}_1(m^i,t))$ converges in distribution
to
$$
\left(\bar{Z}_0(0)+\lambda_0t - \mu_0 \int_0^t \bar{\phi}_0\left(\frac{\bar{Z}_0(s)}{\bar{Z}_0(s)+\bar{Z}_1(s)}\right)\diff s,\
\bar{Z}_1(0)+\lambda_1t - \mu_1 \int_0^t \frac{\bar{Z}_1(s)}{\bar{Z}_0(s)+\bar{Z}_1(s)}\diff s
\right)
$$
and we conclude that, almost surely,
\begin{align*}
    \bar{Z}_0(t) &= \bar{Z}_0(0) + \lambda_0 t - \mu_0 \int_0^t
    \bar{\phi}_0 \left(\frac{\bar{Z}_0(s)}{\bar{Z}_1(s)+\bar{Z}_0(s)}\right)\diff s,\\
    \bar{Z}_1(t) &= \bar{Z}_1(0) + \lambda_1 t - \mu_0 \int_0^t
    \frac{\bar{Z}_1(s)}{\bar{Z}_1(s)+\bar{Z}_0(s)}\diff s,\\
     \bar{Z}_k(t) &= 0,\quad \text{for } 2\leq k \leq L
\end{align*}
holds for all $t$ in $\R_+$.

\begin{lemma}
 $$
\mathcal{C} = \left \{ (\bar{N}(m^i,t)),\ i\in\N \right \}
$$
 is relatively compact and its limiting points are continuous
 processes.
 \label{lemma:relat-compact}
\end{lemma}
\bp The method used here is also standard.
We define $w_h$ the modulus of continuity for any function $h$
defined on $[0,T]$:
$$
w_h(\delta) = \sup_{s,t\leq T;\ |t-s|<\delta} |h(s)-h(t)|.
$$

As in the proof of Theorem \ref{theo:large-numbers}, we can prove that
$(\bar{M}(m,t))$ is a martingale and its increasing processes are given by,
for $0\leq k, l \leq K$ with $k\neq l$,
\begin{align*}
\langle \bar{M}_k(m,t) \rangle &= \frac{\lambda_k t}{\|m\|} + \frac{\mu_k t}{\|m\|}
\int_0^t \phi_k(\bar{N}_0(m,s),\bar{N}_1(m,s),N_2(\|m\|s),\dots,N_L(\|m\|s))\diff s,\\
\langle \bar{M}_k(m,t), \bar{M}_l(m,t)  \rangle &= 0.
\end{align*}
Since the functions $\phi_k$ are bounded by 1, we deduce that
$$
\langle \bar{M}_k(m,t) \rangle \leq \frac{(\lambda_k+\mu_k)t}{\|m\|}.
$$
By using Doob's inequality, we have that, for $\varepsilon>0$,
\begin{align*}
\Pr\left( \sup_{0\leq s \leq t} \|\bar{M}(m,s)\| \geq \varepsilon \right)
&\leq \sum_{k=0}^L \Pr\left( \sup_{0\leq s \leq t} |\bar{M}_k(m,s)| \geq \frac{\varepsilon}{L} \right),\\
&\leq \frac{L^2t \sum_{k=0}^L\lambda_k+\mu_k}{\varepsilon^2 \|m\|}
\end{align*}
We deduce that if $(m^i)$ is a sequence satisfying Equation \eqref{eq:mi}
  then $(\bar{M}(m^i,t))$ converges in probability to 0 uniformly on
  compact sets when $i$ tends to infinity; for any $T\leq 0$ and any
  $\varepsilon>0$,
$$
\lim_{i\rightarrow\infty}\Pr\left(\sup_{0\leq t \leq T}
  \|\bar{M}(m^i,t)\| \geq \varepsilon \right) =0.
$$

Using Equation \eqref{eq:mbar} and  the previous equation, we can prove that
for any $\varepsilon>0$ and $\eta>0$, there exist $\delta>0$ and $A$ such that
for all $i\geq A$,
$$
\Pr\left( w_{\bar{N}_k(m^i,.)}(\delta) > \eta \right) \leq \varepsilon.
$$

The conditions of \cite[7.2 p81]{billingsley-99} are then fulfilled
and the set $\mathcal{C}$ is relatively compact and its limiting points are
continuous processes.
\ep

We are now able to characterize $(\bar{Z}_2,\dots,\bar{Z}_L(t))$.

\begin{lemma}
We consider a sequence $(m^i)$ satisfying Equation \eqref{eq:mi} and such that
$(N(m^i,t))$ converges in distribution to $(\bar{Z}(t))$.
Then the process $(\bar{Z}(t))$ verifies $\bar{Z}_k(t) = 0$ for $2\leq k \leq L$ and $t\geq 0$.
\label{lemma:2-K}
\end{lemma}
\bp
Define the process $(\hat{N}_{2:L}(t))$ with the following transition
rates for $2\leq k \leq L$:
\begin{align*}
  n_k \mapsto n_k+1 &: \lambda_k,\\
  n_k \mapsto n_k-1 &: \mu_k \frac{n_ka_k}{1+n_ka_k}.
\end{align*}
and such that $\hat{N}_{2:L}(0)\geq\tilde{N}^\alpha_{2:L}(0)$.
Clearly, $(\hat{N}_{2:L}(t))$ is ergodic and stochastically dominates $(N_{2:L}(t)$.
Moreover, the components evolve independently.
If we call $(\hat{Z}_{2:L}(t))$ a fluid limit of $(\hat{N}_{2:L}(t))$, we can prove, as in \cite[Prop 5.16 p125]{robert-03},
that it satisfies
$$
\hat{Z}_k(t) = (\hat{Z}_k(0)+(\lambda_k-\mu_k)t)_+,\quad \text{for } t\geq0\text{ and } 2\leq k\leq L.
$$
Moreover, by stochastic domination, we have $\bar{Z}_k(t) \leq \hat{Z}_k(t)$ for
all $t\geq 0$ and $2\leq k \leq L$.
Since $(m^i)$ satisfies Equation \eqref{eq:mi}, we have that $\bar{Z}_k(0)=0$,
for $2\leq k \leq L$ and we conclude that $\bar{Z}_k(t)=0$ for all $t\geq 0$
and $2\leq k \leq L$.
\ep

We define $\bar{\N} = \N \cup \{ +\infty \}$ and we consider
$\bar{\N}^{K-2}$. We endow $\S$ with the metric induced from the $L_1$-norm on
$\R^{K-2}$ by the mapping $(x_2,\dots,x_L)\mapsto (1/(x_2+1),\dots,1/(x_L+1))$.
We can note that, in particular, $\S$ is compact. In the same way as
in \cite{kurtz-94}, we then define a family of random measures on
$[0,\infty)\times\S$, for any $m\in\N^2$, any $\Gamma \subset \S$ and $t\geq0$,
$$
\nu_{m}((0,t)\times\Gamma) = \int_0^t \ind_{\{(N_{2}(\|m\|s),\dots,N_L(\|m\|s))
  \in \Gamma\}}\diff s.
$$
We then consider $\mathcal{L}_0(\S)$, the set of measures $\gamma$
defined on $[0,t]\times\S$ and such that $\gamma([0,t)\times \S) = t$
for all $t\geq0$.  Under the topology induced by weak convergence on
every compact and since $\S$ is compact, $\mathcal{L}_0(\S)$ is
compact. The next lemma shows that a random measure in
$\mathcal{L}_0(\S)$ can be expressed as the sum of probability
measures of $\S$ indexed by $s$ in $\R_+$.

\begin{lemma}
 We consider a sequence $(m^i,\ i\in\N)$ satisfying Equation \eqref{eq:mi}.
The set
$$
\left\{ ((\bar{N}(m^i,t)),\nu_{m^i}), i\in\N \right\}
$$
is relatively compact.  If $((\bar{Z}(t)),\nu)$ is a limit process
then there exists a process $\vartheta$
such that for all $t$, $\vartheta(t,.)$ is a
random probability measure on $\S$ and
$$
\forall t\geq0,\ \forall \Gamma \subset \S,\ \nu([0,t)\times\Gamma)
=\int_0^t \vartheta(s,\Gamma)\diff s.
$$

\label{lemma:measure-lim}
\end{lemma}
\bp One can find a related result in a slightly different context in
\cite{kurtz-92} and \cite{kurtz-94}.

In order to prove the relative compactness of $\left\{
  ((\bar{N}(m^i,t)),\nu_{m^i})\right\}$, we just need to prove that $\left\{
  (\bar{N}(m^i,t)), i\in\N\right\}$ and $\left\{ \nu_{m^i},
  i\in\N\right\}$ are relatively compact. We have already proved the
relative compactness of the first one in lemma
\ref{lemma:relat-compact}. $\mathcal{L}_0(\S)$ is compact then the
second one is relatively compact.

We consider a convergent sequence $((\bar{N}(m^i,t)),\nu_{m^i})$ and
its limit process $((\bar{Z}(t),\nu)$.  Let $(\Omega,\mathcal{F},\Pr)$
be the probability space on which they are defined. We call $\{
\mathcal{F}_t\}$ the natural filtration of $((\bar{Z}(t),\nu)$.

We then define $\gamma$ such that
$$
\forall B \in \mathcal{F},\ \forall C \in
\mathcal{B}([0,\infty))\otimes\mathcal{B}(\S)\ \gamma(B\times C) =
\E(\ind_B\nu(C)).
$$

According to \cite[appendix 8]{ethier-86}, $\gamma$ can be extended to a
measure on
$\mathcal{F}\otimes\mathcal{B}([0,\infty))\otimes\mathcal{B}(\S)$ and there exists
$\vartheta$ such that for all $t$, $\vartheta(t,.)$ is a random probability
measure on $\S$
and for any $B\in \mathcal{B}(\S)$, $(\vartheta(t,B), t\geq 0)$ is
$\{\mathcal{F}_t\}$-adapted and for any $C\in\mathcal{F}\otimes\mathcal{B}([0,\infty))$,
$$
  \gamma(C\times B)
  =\E \left( \int_0^{+\infty} \ind_C(s)\vartheta(s,B)\diff s
  \right).
$$

We define
$$
M_B(t) = \nu([0,t]\times B) - \int_0^t \vartheta(s,B)\diff s.
$$
$((M_B(t))$ is $\{\mathcal{F}_t\}$-adapted and continuous.  We
consider $D\in \mathcal{F}_t$. We define
$\ind_C(\omega,s)=\ind_D(\omega)\ind_{[t,+\infty)}(s)$ and we have
\begin{align*}
  \E\left( \ind_D \nu([t,+\infty)\times B) \right) &= \gamma(C\times B),\\
  &= \E\left( \ind_D \int_t^{+\infty} \vartheta(s,B)\diff s\right).
\end{align*}
It follows that
$$
\E\left( \nu([t,+\infty)\times B) | \mathcal{F}_t \right) = \E\left(
  \int_t^{+\infty} \vartheta(s,B)\diff s | \mathcal{F}_t \right).
$$
Then, $(M_B(t))$ is a continuous $\{\mathcal{F}_t\}$-martingale. It
has finite sample paths and then is almost surely identically
null. Almost surely, the following equation holds for all $t$,
$$
\forall B\subset \S,\ \nu([0,t)\times B) = \int_0^t \vartheta(s,B)\diff s.
$$
\ep

The previous lemma gives a convenient way to express the integral of
linear combinations of indicator functions against a limit measure of
$(\nu_m)$. We will use this technical lemma in the next one to characterize the
limit of $(\nu_m)$. This is the most important part of the theorem and we
exhibit the time scale separation between classes 0 and 1 and $2,\dots,L$ here.

\begin{proposition}
  We consider a sequence $(m^i,\ i\in\N)$ satisfying Equation \eqref{eq:mi} and such that
  $((\bar{N}(m^i,t)),\nu_{m^i})$ is a converging sequence
  and $h$ a continuous function on $[0,1]^2\times\S$. Then, we have that
$$
\left( \int_0^t h(\bar{N}_0(m^i,s),\bar{N}_1(m^i,s),N_2(\|m^i\|s),\dots,N_L(\|m^i\|s))\diff s\right)
$$
converges in distribution to
$$
\left( \int_0^t \sum_{y \in \S}
  h\left(\bar{Z}_0(s),\bar{Z}_1(s),y\right)
  \tilde{\pi}^{\alpha(s)}(y)\diff s \right)
$$
with
$$
\forall t \in \R_+,\ \alpha(t) = \frac{\bar{Z}_0(t)}{\bar{Z}_0(t)+\bar{Z}_1(t)}
$$
and where, $\tilde{\pi}^{\alpha(t)}$ is the
stationary distribution of the process $(\tilde{N}_{2:L}^{\alpha(t)}(s))$
defined in Section \ref{subsection:frozen}.

In particular, the function
$$(x_0,x_1,n_2,\dots,n_L)\mapsto\psi_0(x_0,x_1,n_2a_2,\dots,n_La_L)$$
is continuous on $[0,1]^2\times\S$ and
$$
\left( \int_0^t \phi_0(\bar{N}_0(m^i,s),\bar{N}_1(m^i,s),N_2(\|m^i\|s),\dots,N_L(\|m^i\|s))\diff s\right)
$$
converges in distribution to
$$
\left( \int_0^t \bar{\phi}_0\left(\alpha(s)\right)\diff s \right)
$$
\label{prop:averaging}
\end{proposition}
\bp
The functions which are continuous on $[0,1]^2\times\S$ for the topology induced
by the natural topology on $[0,1]^2$ and the topology induced by the mapping
$(x_2,\dots,x_L)\mapsto(1/(x_2+1),\dots,1/(x_L+1))$ are the bounded functions $h$
which are continuous on $[0,1]^2\times\N^{K-2}$ for the natural topology and such that
for each $x\in[0,1]^2$, $y \mapsto h(x,y)$ admits a unique limit $h(x,\infty)$
in all directions such that $\|y\|\to\infty$ and the function $x\mapsto h(x,\infty)$
has to be continuous on $[0,1]^2$. We can remark that, for all $(x_0,x_1)\in[0,1]^2$
$\psi_0(x_0,x_1,n_{2:L}\hadam a_{2:L})\to0$ when $\|n_{2:L}\|\to+\infty$ and the function
$$(x_0,x_1,n_2,\dots,n_L)\mapsto\psi_0(x_0,x_1,n_2a_2,\dots,n_La_L)$$
is then continuous on $[0,1]^2\times\S$.

We consider $h$ a continuous function on $[0,1]^2\times\S$, since the space
$[0,1]^2\times\S$ is compact, Lemma \ref{lemma:measure-lim} implies directly that
$$
\left( \int_0^t h(\bar{N}_0(m^i,s),\bar{N}_1(m^i,s),N_2(\|m^i\|s),\dots,N_L(\|m^i\|s))\diff s\right)
$$
converges in distribution to
$$
\left( \int_0^t \sum_{y \in \S}
  h\left(\bar{Z}_0(s),\bar{Z}_1(s),y\right)
  \vartheta(s,y)\diff s \right).
$$

We are now able to fully characterize the random measures $(\vartheta(.,t))$.
For any continuous bounded function $f$ on $\S$ and any $m\in\N^{K-2}$, we define,
\begin{align*}
  \bar{M}_f(m,t) = \frac{1}{\|m\|} \Bigl( f(N_2&(\|m\|t),\dots,N_L(\|m\|t)) - f(0)\Bigr)\\
  - \sum_{k=2}^L\lambda_k \int_0^t &\Bigl(f((N_2(\|m\|s),\dots,N_L(\|m\|s))+e_k) -
  f(N_2(\|m\|s),\dots,N_L(\|m\|s)) \Bigr) \diff s  \\
   - \sum_{k=2}^L \mu_k \int_0^t &\Bigl(
      \bigl(f((N_2(\|m\|s),\dots,N_L(\|m\|s))-e_k) - f(N_2(\|m\|s),\dots,N_L(\|m\|s)) \bigr)\\
    &\phi_k(\bar{N}_0(m,s),\bar{N}_1(m,s),N_2(\|m\|s),\dots,N_L(\|m\|s))
  \Bigr) \diff s.
\end{align*}

As $(\bar{M}(m,t))$ defined by equation (\ref{eq:mbar}) is a
martingale, $(\bar{M}_f(m,t))$ is a martingale. We consider a
convergent sequence $((\bar{N}(m^i,t)),\nu_{m^i})$. We have that
$(\bar{M}_f(m^i,t))$ converges in distribution to $0$. $\|m^i\|^{-1}
(f(N_2(\|m^i\|t),\dots,N_L(\|m^i\|t)) - f(0))$ also converges to $0$ because
$f$ is bounded. As a consequence, the following term
\begin{align*}
  \sum_{k=2}^L \lambda_k \int_0^t
  \Bigl(&f((N_2(\|m^i\|s),\dots,N_L(\|m^i\|s))+e_k) -
    f(N_2(\|m^i\|s),\dots,N_L(\|m^i\|s)) \Bigr)\diff s  \\
  -\sum_{k=2}^L \mu_k \int_0^t \Bigl(& \left
      (f((N_2(\|m^i\|s),\dots,N_L(\|m^i\|s))-e_k) - f(N_2(\|m^i\|s),\dots,N_L(\|m^i\|s))
    \right)\\
    &\phi_k(\bar{N}_0(m^i,s),\bar{N}_1(m^i,s),N_2(\|m^i\|s),\dots,N_L(\|m^i\|s))\Bigr)\diff s
\end{align*}
also converges in distribution to $0$. But, by the continuous mapping theorem and Lemma
\ref{lemma:measure-lim}, it converges in distribution to
\begin{align*}
\int_0^t \sum_{k=2}^L \Biggl ( &\lambda_k \sum_{y \in \S} f(y+e_k) -
  f(y)\\ &+ \mu_k \sum_{y \in \S} (f(y-e_k) - f(y))
  \varphi_k(\bar{Z}_0(s),\bar{Z}_1(s),y)\Biggr) \vartheta(s,y)\diff s.
\end{align*}
Consequently, this is null almost surely for all $t$ and we have then, for
Lebesgue-almost every $t$,
\begin{align*}
\sum_{k=2}^L \Biggl ( &\lambda_k \sum_{y \in \S} f(y+e_k) - f(y) +\\
 &\mu_k \sum_{y \in \S} (f(y-e_k) - f(y))
  \varphi_k(\bar{Z}_0(t),\bar{Z}_1(t),y)\Biggr) \vartheta(t,y) =0.
\end{align*}
We deduce immediately that
$$
\int_{\N^{K-2}}
\tilde{\Omega}^{\alpha(t)}(f)(y)\vartheta(t,\diff y)=0
$$
where $\alpha(t)=\bar{Z}_0(t)/(\bar{Z}_0(t)+\bar{Z}_1(t))$ and
$\tilde{\Omega}^{\alpha(t)}$ is the infinitesimal generator of
$(\tilde{N}_{2:L}^{\alpha(t)}(s))$, defined in Section \ref{subsection:frozen}.
This proves exactly that $\vartheta(t,.)$ is invariant for
$(\tilde{N}_{2:L}^{\alpha(t)}(s))$.
By uniqueness of the invariant distribution of such a process, we have that
$$
\vartheta(t,.)=\tilde{\pi}^{\alpha(t)}.
$$
\ep